\documentclass[10pt]{article}
\usepackage{amssymb}
\usepackage{latexsym}
\usepackage{euscript}
\usepackage[frame,cmtip,arrow,matrix,line,graph]{xy}

\newcommand{\pp}{\mathbb{P}}
\newcommand{\qq}{\mathbb{Q}}
\newcommand{\cc}{\mathbb{C}}
\newcommand{\rr}{\mathbb{R}}

\newcommand{\zz}{\mathbb{Z}}
\newcommand{\sss}{\mathbb{S}}

\DeclareFontFamily{OT1}{rsfs}{}
\DeclareFontShape{OT1}{rsfs}{n}{it}{<-> rsfs10}{}
\DeclareMathAlphabet{\curly}{OT1}{rsfs}{n}{it}

\newcommand{\cd}{\partial}
\newcommand{\cdbar}{\overline{\partial}}

\newcommand{\im}{\mathrm{im \,}}

\newcommand{\coker}{\mathrm{coker \,}}

\newcommand{\Sym}{\mathrm{Sym}}

\newcommand{\Crit}{\mathrm{Crit}}

\newcommand{\Pic}{\mathrm{Pic}}
\newcommand{\ind}{\mathrm{Index \,}}
\newcommand{\id}{\mathrm{id \,}}

\newcommand{\Fibre}{\mathrm{Fibre}}

\newcommand{\Jac}{\mathrm{Jac}}
\newcommand{\Hom}{\mathrm{Hom}}
\newcommand{\rk}{\mathrm{rk}}

\newcommand{\mgbar}{\overline{M}_g}

\newcommand{\oz}{\overline{z}}
\newcommand{\usigma}{\underline{\sigma}}
\newcommand{\ueta}{\underline{\eta}}
\newcommand{\cok}{\mathrm{cok}}
\newcommand{\Hilb}{\mathrm{Hilb}}
\newcommand{\norm}{\mathrm{norm}}

\newcommand{\LL}{\mathcal{L}}

\newcounter{Universal}[section]
\renewcommand{\theUniversal}{\thesection.\arabic{Universal}}

\newenvironment{Plain}{\refstepcounter{Universal} \par \vspace{0.5cm}
\noindent {\bf (\theUniversal)}\ }{\par \vspace{0.5cm}}
\newenvironment{Italic}{\refstepcounter{Universal} \par \vspace{0.5cm}
\noindent {\bf (\theUniversal)}\ \it}{\par \vspace{0.5cm}}  
  
\newenvironment{Thm}{\begin{Italic}{\sc Theorem: }}{\end{Italic}}
\newenvironment{Prop}{\begin{Italic}{\sc Proposition:} }{\end{Italic}}
\newenvironment{Defn}{\begin{Italic}{\sc Definition: } }{\end{Italic}}  
\newenvironment{Cor}{\begin{Italic}{\sc Corollary: }}{\end{Italic}}
\newenvironment{Lem}{\begin{Italic}{\sc Lemma: }}{\end{Italic}}

\newenvironment{Pf}{\par \noindent{\sc Proof:} }{\quad $\blacksquare$ \par
\vspace{0.5cm}}
\newenvironment{Example}{\begin{Plain} \noindent{\sc Example: }}
{\quad $\square$ \end{Plain}} 
\newenvironment{Rmk}{\begin{Plain} \noindent{\sc Remark: }}{\quad
    $\square$ \end{Plain}}

\newenvironment{Question}{\begin{Plain} \noindent{\sc Question:
}}{\quad $\square$ \end{Plain}}

\newcounter{enum}

\newenvironment{Eqn}{\refstepcounter{Universal} $$} {\eqno \mathrm{
(\theUniversal)} $$} 

\title{Lefschetz pencils and the canonical class for symplectic $4$-manifolds}
\author{Simon Donaldson\thanks{Imperial College, London: s.donaldson@ic.ac.uk} \ and Ivan
  Smith\thanks{New College, Oxford: smithi@maths.ox.ac.uk}}
\date{}

\begin{document}
\maketitle
\thispagestyle{empty}

\begin{abstract}
\noindent We present a new proof of a result due to Taubes:  if $(X,\omega)$ is
a closed  symplectic four-manifold with $b_+ (X) > 1+b_1 (X)$ and
$\lambda[\omega] \in H^2 (X; \qq)$ for some $\lambda \in \rr^+$, then the
Poincar{\'e} dual of $K_X$ may be represented by an embedded symplectic
submanifold.  The result builds on the existence of
Lefschetz pencils on symplectic four-manifolds.  We approach the topological
problem of constructing submanifolds with locally positive
intersections via almost-complex geometry.  The crux of the argument
is that a Gromov invariant counting pseudoholomorphic sections of an
associated bundle of symmetric products is non-zero. 
\end{abstract}


\section{Introduction}
 
Taubes' renowned work on the Seiberg-Witten invariants of 
symplectic $4$-manifolds shone light on an area which had, before then, been 
largely the domain of conjecture and speculation. The conclusions of Taubes' 
work can be divided into two categories. First, results relating the existence 
of a symplectic structure to the differential topology of the underlying 
$4$-manifold; second, results which stay within the framework of symplectic 
geometry. The results in the first category hinge on the fact 
that the Seiberg-Witten theory yields invariants of the smooth structure, and 
the Seiberg-Witten theory, or something equivalent, seems unavoidable. For the 
results in the second category, on the other hand, it is natural to ask if 
there are alternative proofs, avoiding the Seiberg-Witten theory and
staying more  
within the realm of symplectic geometry. The purpose of this paper is
to take a step in this direction. Specifically, we will give a new
proof of the

\begin{Thm} [Taubes] \label{maintheorem}
Let $(X, \omega)$ be a closed symplectic four-manifold with $b_+ (X) > 1+ b_1
(X)$ and such that some positive real multiple of the cohomology class
$[\omega] \in H^2(X; \rr)$ is
rational. Then there is a smooth embedded symplectic surface in $X$
which represents the homology class Poincar{\'e} dual to $K_X$.
\end{Thm}

\noindent Here, as usual, $b_{+}(X)$ is the dimension of a maximal
positive subspace for  
the intersection form on $X$, endowed with the conventional orientation, 
while $b_{1}$ is the first Betti number. 
The 
symplectic structure defines, up to homotopy, an almost-complex structure on 
$X$, and $K_{X}\in H^{2}(X;\zz)$ is the 
first Chern class of the cotangent bundle. A surface $\Sigma\subset X$ is 
\emph{symplectic} if the restriction of $\omega$  to $\Sigma$ is
nondegenerate. The statement  
of (\ref{maintheorem}) is weaker than the result proved by Taubes
(\cite{Taubes}, \cite{Taubes:SWtoGR}, \cite{Taubesproof}) in two respects.  
First, Taubes' proof only requires the hypothesis that $b_{+}>1$.  As
we explain later, the
tighter constraint simplifies the proof but does not seem essential,
and one can hope 
to obtain the stronger statement by a suitable refinement of the techniques 
we develop in this paper. Second, we impose the requirement that (a
multiple of) $[\omega]$ be 
a rational class, which is not required by Taubes.
It is not possible to remove this without some
significant extension 
of the arguments, as we explain in more detail below (and see
\cite{Sequel});  however, the
restriction does not seem significant for most applications, since
every symplectic form is deformation equivalent to a
rational symplectic form.

\vspace{0.2cm}

\noindent Our proof of (\ref{maintheorem}) will avoid any use of the
Seiberg-Witten equations - although see the remarks below - and is in
general ``softer'' (in 
Gromov's terminology \cite{GromovIC}) than
Taubes'. Thus we avoid 
delicate arguments with 
elliptic partial differential equations, although elliptic theory does 
play  a vital role in 
our argument: at one point through the index theory for linear operators and 
at another point through pseudo-holomorphic curves and  a  ``Gromov 
invariant''. 

\vspace{0.2cm}

\noindent Symplectic $4$-manifolds occupy the ground between, on one
side, general  
smooth $4$-manifolds and, on the other side, complex algebraic surfaces.  In 
this spirit, one path by which to approach Taubes' result is to recall the 
proof in the case when $X$ is a complex algebaic surface and $\omega$ is a 
K{\"a}hler form. Then one can simply argue as follows. The Hodge index formula 
states that $b_{+}= 1 + 2 p_{g}$ so the hypothesis $b_{+}>1$ says $p_{g}>0$, 
i.e. there is a non-trivial holomorphic $2$-form $\Theta$ on $X$. The zero-set 
of $\Theta$ is a complex curve representing $K_{X}$ and, at least if this 
curve is non-singular, this immediately gives the result. (In the rare event 
 that, for all $\Theta$, the curve is 
singular, one would need to consider a $C^{\infty}$ perturbation, in
the manner  
of Proposition (\ref{smoothinglemma}) below.)  Our proof - and, to some
extent, that of Taubes - can be seen as 
a translation of this simple argument into the symplectic setting. In this 
case we  always have compatible \emph{almost-complex} structures, and the crux
of the matter is to get around the lack of integrability.
 
\vspace{0.2cm}

\noindent The starting point for our proof is the existence, proved in
\cite{Donaldson:pencils}, of Lefschetz  
pencils on symplectic $4$-manifolds. This means that a blow-up $X'$ of 
$X$ is the total space of a Lefschetz fibration $\pi: X' \rightarrow 
S^{2}$. Here the fibres are Riemann surfaces, so both the base and fibre are 
complex and the deviation from integrability is confined to the twisting of 
the fibration. In the classical case, with a complex surface and holomorphic 
Lefschetz fibration, standard techniques allow us to study the geometry of the 
total space, and in particular the holomorphic $2$-forms, in terms of data on 
the base and fibres---in algebro-geometric language we consider direct image 
sheaves over the Riemann sphere. Our proof can be understood, from one point 
of view, as an adaptation of these arguments to the non-integrable case, and 
we obtain our symplectic surface via suitable sections of a bundle
$X_r(f)$, with $r = 2g-2$, of 
symmetric products along the fibres. 

\vspace{0.2cm}

\noindent A second point of view is to try to trace the parallel
between our arguments  
and those of Taubes, using the Seiberg-Witten equations. Thus we ask
what might be  
said about the solutions of the Seiberg-Witten equations over the total space 
of a Lefschetz fibration. Here we make contact with a programme which is being 
developed by D. Salamon and his collaborators involving the  
``adiabatic limit'' of the Seiberg-Witten equations. In that programme
one would  
consider a family of metrics $g_{\epsilon}$ on $X'$ in which the fibre 
has diameter $O(\epsilon)$, and study the asymptotics as $\epsilon\rightarrow 
0$. One would expect, by analogy with earlier work of Dostoglou and
Salamon \cite{Dost-Sal}, that when $\epsilon$ is small the
Seiberg-Witten solutions are modelled on  
holomorphic sections of an associated bundle of ``vortices'' on the 
fibres (the 
solutions of the Seiberg-Witten equations on a Riemann surface times 
$\rr^{2}$, invariant under translations in the $\rr^{2}$ variables). 
It is a well-known and simple result that the moduli spaces of 
vortices can be identified with  the symmetric products of 
the Riemann surface, so one would expect to arrive
at a similar picture to that studied in 
the present paper.  (Indeed in the introduction to $Gr \Rightarrow
SW$, which appears as the third paper in \cite{Taubesproof},
Taubes introduces an elliptic operator which looks at vortices in the
normal bundle fibres to a fixed holomorphic curve, and which he describes as
interpolating between the Gromov and Seiberg-Witten theories.  Related
work has also been done by Hong, Jost and Struwe \cite{Jostetal}.)
Put another way, the almost-complex structure on
$X'$ defines a natural almost-complex structure on $X_r(f)$:  the tangent
space to $X_r(f)$ at a tuple $(p_1 + \cdots
+ p_r)$ is an extension of the direct sum $\oplus T_{p_i} \Sigma$ by a
tangent direction to the base, and each component
has a given anti-involution $J$.  This almost-complex
structure is $C^0$ but not $C^1$;
nonetheless, it is the limit of smooth structures (this is exploited
in \cite{Sequel}).  The
parameter $\epsilon$ plays a similar 
role to the deformation parameter used by Taubes, and in both cases the main 
difficulties have to do with the details of the asymptotic analysis with 
respect to the parameter. 

\vspace{0.2cm}

\noindent The advantage of our argument, in which the 
connection with the Seiberg-Witten theory is kept in the background, is that 
we do not need to contend with such problems. From this perspective,
the simplification has two fundamental aspects.   Just as Taubes, we regard
a surface in a four-manifold as the zero set of equations $f(z,w) =
0$.  Since the situation is not integrable, if $z$ and $w$ are local
almost-complex co-ordinates, then 
we cannot expect $f$ to be holomorphic in both variables.  However, by
first imposing a Lefschetz pencil, we
can ensure that it is \emph{holomorphic} in $z$ (along the fibres of a
pencil) and smooth in $w$ (varying the element of the pencil).
Second, we replace the \emph{non-compact} deformation of a metric in Taubes
proof with a \emph{compact} deformation of almost-complex structures, in the
definition of the Gromov invariant.  Some additional consideration of the
degenerating family of smooth almost complex structures mentioned
above, or the 
Salamon programme, would be needed to remove the hypothesis that the
symplectic form be rational from our approach.

\vspace{0.2cm}

\noindent  Taubes' proof of (\ref{maintheorem})
can be seen as falling into two parts.
One part is the identification of the Seiberg-Witten invariant of a line 
bundle $L\rightarrow X$ with 
the Gromov invariant counting pseudo-holomorphic curves homologous to
$(K_{X}/2) + c_{1} (L)$ (assuming for simplicity that $K_{X}$ 
is even). This uses Taubes' perturbation of the Seiberg-Witten 
equations. The second part is the statement that the invariants for $L$ and 
$L^{*}$ are equal. This symmetry is completely trivial using the original, 
unperturbed, Seiberg-Witten equations, but in conjunction with the first part 
gives a highly non-obvious assertion about pseudo-holomorphic curves:

\vspace{0.2cm}

\noindent \emph{The Gromov invariants defined by counting holomorphic curves in
the homology classes $(K_{X}/2) + \lambda$ and $(K_{X}/2) -\lambda$
are equal.}  

\vspace{0.2cm}
 
\noindent Taking $\lambda= (K_X /2)$, one of the pair
of homology 
classes is zero which trivially - or by definition - has Gromov invariant $1$
(the empty set is the unique pseudo-holomorphic representative).  Hence
this symmetry includes as a special case Taubes' result that the Gromov 
invariant of $K_{X}$ is $1$, which (again 
up to the issue of perhaps smoothing singular pseudo-holomorphic curves) gives 
(\ref{maintheorem}).  In our setting the fundamental symmetry appears
through the duality between divisors $D$ and $K_C - D$ on a Riemann
surface $C$.  A detailed treatment of this, and of some constructions of
symplectic surfaces in four-manifolds with $b_+ = 1$, can be found in
the paper \cite{Sequel}.

\vspace{0.3cm}

\noindent \textbf{Outline of the paper:}

\begin{enumerate}

\item As indicated above, we will produce the desired symplectic
  surface from a section of a (compactified) bundle of symmetric
  products.  The proof hinges on the availability of two almost
  complex structures on this space - one which allows the computation
  of a Gromov invariant, and one for which the resulting
  pseudoholomorphic sections will yield embedded symplectic surfaces.

\item  After a quick review of Lefschetz pencils, we introduce
  ``positive symplectic divisors'' and prove a smoothing lemma which
  gives sufficient conditions to deform such to embedded symplectic
  surfaces (section 2).

\item Given a Lefschetz fibration $f: X \rightarrow \sss^2$ there is an
  associated fibration 
  $X_r(f)$ whose fibre over a smooth point $t \in \sss^2$ is the
  $r$-th symmetric 
  product of the fibre of the given pencil. This bundle of symmetric
  products can be naturally compactified to a smooth 
  space (section 3 and the Appendix).   

\item Write $2g-2=r$ henceforth. The fibres of the Lefschetz pencil
  are Riemann surfaces.  The vector spaces $H^0 (\Sigma_t;
  K_{\Sigma_t})$ patch to give a vector bundle $V$ over $\sss^2$ whose
  projectivisation naturally embeds into $X_r(f)$,
  avoiding the locus of 
  critical values of $X_r(f) \rightarrow \sss^2$ (section 3 and the Appendix).

\item A nowhere zero section of $V$ defines a two-dimensional homology
  class $[\psi_V]$ in
  $X_r(f)$ via the embedding of $\pp (V)$.  $X_r(f)$ admits symplectic
  structures and there is a well-defined
  Gromov invariant which counts pseudoholomorphic sections of $X_r(f)$ in
  the class $[\psi_V]$, a problem of virtual index zero (section 4).

\item Fixing
  a generic $\cdbar$-operator on $V$ gives a distinguished
  $N$-dimensional space of 
  holomorphic sections, for $N =  [b_+(X) -1 - b_1(X)]/2$ (section
  5). For a suitable almost-complex structure, the moduli space of
  holomorphic sections of $X_r(f)$ in the distinguished homology class is
  the projective space $\pp^{N-1}$.  To prove this one considers a
  natural map from $X_r(f)$ to a fibration of Picard varieties $P_r (f)$,
  for which sections are generically precluded by the
  existence of sections of $V$ (section 5).

\item The Gromov invariant above can be
  evaluated as the Euler class of an
  obstruction bundle over the projective space; the obstruction bundle
  is the dual 
  of the quotient bundle, and the Gromov invariant is $\pm 1$ 
  (section 5).

\item The symmetric product fibration contains a number of natural
  diagonal strata; although singular, these have ``smooth models'',
  and there are almost-complex structures on $X_r (f)$ which are
  compatible with the natural inclusions of the strata (section 6).

\item For such an almost-complex structure on $X_r(f)$, a
  pseudoholomorphic section $\phi$ gives rise to a positive symplectic
  divisor $C_{\phi}$ inside $X'$. 
  If the original Lefschetz pencil is of sufficiently high degree,
  $\phi$  has no bubble components, and the divisor $C_{\phi}$
  descends to $X$, where it may be smoothed to give an embedded
  symplectic surface in the class $PD[K_X]$ (section 7).
\end{enumerate}

\noindent It may be worth stressing that the arguments of the paper go
through with comparatively little machinery from the
pseudoholomorphic curves industry.  Whilst the ``virtual class''
methods have brought the theory of
Gromov-Witten invariants to a very 
satisfactory status, the more direct arguments of \cite{McD-S:Jhol} or
Gromov's original \cite{Gromov} suffice in this paper. 

\vspace{0.2cm}

\noindent \textbf{Acknowledgements:} The second author is grateful to
Paul Seidel and Richard Thomas for helpful conversations.  He is also
indebted to the first author for permission to present evolving aspects of
the work in conferences at the ETH (Zurich: December 1999) and at Ecole
Polytechnique (Paris: June 2000), and thanks the organisers of both
workshops for their hospitality.  Both authors thank the referee for
detailed comments, which have certainly helped improve the paper.


\section{Lefschetz pencils and standard surfaces}

Let $(X, \omega_X)$ be a closed symplectic four-manifold.  As in
\cite{Donaldson:pencils}, a system of local complex co-ordinates $(z_1, z_2)$
centred at a point $x \in X$ is \emph{compatible} with $\omega$ if the
symplectic form at $\underline{z} = 0$ is positive and of type
$(1,1)$.  Local co-ordinates are always assumed compatible with the
given global orientation on $X$.

\begin{Defn} \label{lefpencildefn}
A \emph{topological Lefschetz
pencil} of degree $k$ on $X$ comprises a map $f: X \backslash \{b_1,
\ldots, b_d \}
\rightarrow \sss^2$ defined on the complement of finitely many points
in $X$ to the two-sphere, which is a submersion outside of finitely
many critical points $\{ p_1, \ldots, p_r \}$, all in distinct 
  fibres of $f$.  We demand this data conforms to local models:

\begin{enumerate}
\item at each $p_i$ there are compatible local complex co-ordinates with
  respect to which $f$ has the form $(z_1, z_2) \mapsto z_1 ^2 + z_2
  ^2$;
\item at each $b_i$ there are compatible local complex co-ordinates
  with respect to which $f$ 
  has the form $(z_1, z_2) \mapsto z_1 / z_2$.
\end{enumerate}
\end{Defn}

\noindent   The closures of
the (open) fibres of $f$ are given by including the
points $\{b_1, \ldots, b_d \}$; it follows from the local models that
the fibres are then manifolds near the $b_i \in X$.  The
results of this paper start with a general existence theorem for
Lefschetz pencils \cite{Donaldson:submflds}, \cite{Donaldson:pencils}:

\begin{Thm} \label{pencilsexist} Let $(X,\omega)$ be a
  symplectic four-manifold for which $[\omega/2\pi]
\in H^2 (X; \zz) \subset H^2 (X; \rr)$.  Then for any sufficiently
large $k$ there 
are Lefschetz pencils on $X$ for which the closed fibres are
Poincar{\'e} dual to $k[\omega]$ and are symplectic submanifolds away
from their singularities.  
\end{Thm}

\noindent The adjunction formula implies that the
genus $g$ of a smooth fibre increases with $k$.  In fact the pencils
have an asymptotic uniqueness property, but we shall not make use of
that in this paper. However, the choice of degree does enter into the
argument~:

\begin{Rmk} \label{highdegree}
At various stages of the proof, it will be important that we fix on
the given $(X, \omega)$ a Lefschetz pencil of sufficiently high degree
$k$~; for instance, it would be sufficient to ensure that $g > \max \{
b_+(X)/2, 2 \}$ and $k > 2\pi (\omega \cdot \omega) / (\omega \cdot
K_X)$. We will explain
why these conditions arise at the relevant points.
\end{Rmk}

\noindent  The singular fibres
of a pencil are nodal curves, and \emph{a priori} the node could
separate the fibre.  We note that you can always avoid these reducible
fibres, although the proof does not rely on this in a crucial fashion:

\begin{Prop}
Every symplectic four-manifold admits a Lefschetz pencil in which
every fibre is irreducible.
\end{Prop}

\noindent One route to proving this is by stabilisation, as in
\cite{ivanmodulidivisor}~; it also
follows from the formulae of Auroux and Katzarkov
\cite{AKstabilisation}.  (If $K_X$ or the intersection form of $X$ is
even it is purely elementary.) 
Now fix compatible almost-complex structures on
the tangent spaces to $X$  
at the points $b_{i}$. We can then define the ``blow-up'' $X'$ of $X$ 
at these $d$ points, as a set, in the usual way by replacing the points by the 
complex projectivisations of their tangent spaces. Again in a familiar way we 
can endow this with the structure of a smooth $4$-manifold with a smooth 
blowing-down map
$p:X' \rightarrow X$.  On the other hand, the map $f$ 
clearly induces  a map (which we still call) $f:X' \rightarrow
\sss^{2}$. The exceptional  
spheres  $E_{i}$ in the blow-up appear as  sections of $f$.

\vspace{0.2cm}

\noindent Now for any $N>0$ we can define a $2$-form
$\omega_{(N)}= p^{*}(\omega_X) + N f^{*}(\omega_{st})$ on $X'$,
where $\omega_{st}$ is the standard area form on $\sss^{2}$ with total
area $2\pi$~: that is,  $\int_{\sss^2} \omega_{st} = 2\pi$.

\begin{Lem} \label{symplecticforms}
The form $\omega_{(N)}$ is a symplectic form on $X'$.
There are disjoint embedded balls $B_{i}\subset X$ containing
$b_{i}$ such that the manifold $X' \backslash \bigcup E_{i}$, with the
symplectic  
form $\omega_{(N)}$, is symplectomorphic to $X \backslash \bigcup B_{i}$
with the symplectic form $(1+kN)\omega_X$.
\end{Lem}

\begin{Pf}
Clearly $\omega_{(N)}$ is closed so we need to see it is non-degenerate.
We have 

$$\omega_{(N)}^{2} = p^{*}(\omega_X^{2}) + 2N p^{*}(\omega_X)\wedge 
f^{*}(\omega_{st}).$$

\noindent The first term is non-negative, and strictly positive away
from the $E_{i}$,  
while the second term is strictly positive  away from the points $p_{j}$, so
$\omega_{(N)}^{2}$ is strictly positive everywhere.

\vspace{0.2cm}

\noindent For the second part we identify $X' \backslash \bigcup E_{i}$ with 
$X \backslash \{b_{1},\dots, b_{d}\}$ via the blow down map $p$. Clearly
$H^{2}(X \backslash \{b_{1},\dots ,b_{d}\})$ is isomorphic to $H^{2}(X)$~; 
we claim the de Rham cohomology class represented by the restriction of 
$f^{*}(\omega_{st})$ is $k[\omega_X]$.  It is clear that it is some
linear multiple of $[\omega_X]$, and the constant follows by
integrating over a fibre~:  

$$\int_F \omega_X \ = \ [\omega_X] \cdot [(k/2\pi) \omega_X] \ = \
(2\pi/k) \sharp \{ \mathrm{Basepoints} \}$$

\noindent as compared with

$$\int_F \omega_{st} \ = \ \mathrm{Area}(\sss^2) \cdot \sharp \{
\mathrm{Basepoints} \}.$$ 

\noindent Thus we can write

\begin{Eqn} \label{scaledsymplforms}
\omega_{(N)}= (1+kN) \omega_X + N da 
\end{Eqn}

\noindent on this punctured manifold. The $1$-form $a$ will not extend
over the  
punctures; in fact we may suppose that, in standard complex co-ordinates 
centred on a puncture,

$$a= \frac{1}{4\pi(\vert z_{1}\vert^{2} + \vert z_{2}\vert^{2})}(\oz_{1} 
dz_{1} - z_{1}d\oz_{1} + \oz_{2} dz_{2} - z_{2} d\oz_{2}). 
$$

\noindent Rewrite equation (\ref{scaledsymplforms}) as

$$ \frac{1}{kN+1} \, \omega_{(N)} = \omega_X + t da, $$

\noindent where $t=N/(kN+1)$. If we were working on a compact manifold
we could  
immediately apply Moser's theorem to obtain a symplectomorphism between 
$\frac{1}{kN+1}\omega_{(N)}$ and $\omega_X$, integrating a time-dependent 
vector field. The argument does not apply immediately as things stand, since 
one may not be able to integrate vector fields on a non-compact manifold. 
However closer examination shows that we can carry through the argument since 
the relevant  vector field  points {\it into} the manifold near the
puncture, so we can integrate for {\it positive} time. The basic model is the 
case of $\cc^{2}$ when the map

$$ F_{t}(z_{1},z_{2}) = \sqrt{4\pi t+1/\vert z\vert^{2}} (z_1, z_2)  $$ 

\noindent gives a diffeomorphism from $\cc^{2}\backslash \{0\}$ to the
complement of the 
ball of radius $ \sqrt{4\pi t}$, and

$$ F^{*}_{t}(\Omega_{\cc^{2}}) = \Omega_{\cc^{2}} +t da_{\cc^{2}},$$

\noindent where $a$ is the $1$-form above (cf. Lemma 6.40 in
\cite{McD-S}).
\end{Pf}

\noindent Using this Lemma, we see that to obtain a symplectic surface
in $X$  it suffices to find a surface in $X'$ which does not meet the 
exceptional curves $E_{i}$ and which is symplectic with respect to the form
$\omega_{(N)}$ for some $N$. 

\begin{Defn} \label{standardsurface}
A \emph{standard surface} $\Sigma$ in $X'$ is an embedded,  oriented 
surface disjoint from the critical points $p^{-1}(p_{j})$ and such that 
the restriction $f|_{\Sigma}$ satisfies the following.

\begin{enumerate}
\item The map defines a branched covering from $\Sigma$
to $\sss^{2}$ which has positive degree on each component of $\Sigma$
and which has only simple branch points.
\item At each branch point $q$, $(df)_q$ defines an
  isomorphism from the normal bundle $(\nu_{\Sigma / X'})_q$ to
  $T_{f(q)}\sss^2$ which is orientation-preserving (with respect to
  the orientation induced by the given orientations on $X'$ and
  $\Sigma$).
\end{enumerate}
\end{Defn}

\noindent If $\Sigma$ is a standard surface, with branch points
$\{ q_{j} \}$ the  
restriction of $f$ to $\Sigma \backslash \{q_{j}\}$ is a covering
map, and  
near each branch point $q_{j}$ there are (compatibly oriented)
co-ordinates which identify 
the triple $(X',\Sigma, f)$ with the standard model 

$$ \{ (z,w)\in \cc^{2}: z= w^{2}\} \ ;\ (z,w)\mapsto z\in \cc. $$

\noindent In particular the tangent space of $\Sigma$ at $q_{j}$
is the same as the  
tangent space of the fibre, and the orientations match up. 

\begin{Lem} A standard surface $\Sigma$  in $X'$  is 
     symplectic with respect to the form $\omega_{(N)}$ for large enough $N$.
\end{Lem}

\noindent This is elementary.
The restriction of $p^{*}(\omega_X)$ to $\Sigma$ is positive near the branch 
points, since it is positive on the fibres of $f$, and the restriction of 
$f^{*}(\omega_{st})$ is non-negative everywhere and strictly positive 
away from the branch points.

In sum we see that {\it a standard surface in $X$ which does not meet the 
exceptional curves $E_{i}$ yields a symplectic surface in the original 
manifold $X$.}
There is a converse to this statement which we will state here, although it is 
not directly relevant in the present paper.

\begin{Prop}
Let $\Sigma$ be any
symplectic surface in $X$. Replacing $\omega_X$ by a 
sufficiently high multiple $k\omega_X$, there is a Lefschetz pencil on 
$X$ such that $\Sigma$ arises from a standard surface in the blow-up $X'$, 
by the correspondence above.
\end{Prop}

\noindent  The proof of this Proposition is a straightforward
modification of the proof
of the existence of Lefschetz pencils in \cite{Donaldson:pencils}
(branched covering maps are the 
analogue, in real dimension $2$, of Lefschetz pencils in real dimension $4$).
The upshot of this discussion is that the study of symplectic surfaces in 
symplectic $4$-manifolds is reduced, in principle, to the study of
standard surfaces which can in turn be described in a rather obvious way by
``monodromy'' data.

\vspace{0.2cm}

\noindent In the proof of our main theorem we will combine the
discussion above with a  
simple smoothing criterion. Suppose $\Sigma_{1},\dots,\Sigma_{n}$ are
symplectic  
surfaces in a compact symplectic $4$-manifold $X$ which intersect
transversally with
locally positive intersection numbers (and no triple intersections).

\begin{Defn} \label{possympldivisor}
For positive integers $a_{i}$ we 
call the formal sum $\sum_{i=1}^{n} a_{i} \Sigma_{i}$
a \emph{positive symplectic divisor} in $X$.
\end{Defn}

\noindent The prototypical construction is then~:

\begin{Prop} \label{smoothinglemma}
Suppose $\sum a_{i}\Sigma_{i}$ is a positive symplectic divisor in $X$ and 
that 

$$   \big(\sum_{i} a_{i} [\Sigma_{i}] \big) \cdot [\Sigma_{j}] \geq 0 $$

\noindent for each $j$. Then in an arbitrarily small neighbourhood of
$\bigcup_{i}  
\Sigma_{i}$ there is a smooth symplectic surface representing the homology 
class $\sum_{i} a_{i}[ \Sigma_{i}] $.
\end{Prop} 

\begin{Pf}
Recall, by assumption, that there are no triple intersections.  
Our proof will appeal to a local complex model for the union of the
surfaces, but as things stand such may not exist.  Let $x$ be one of
the transverse intersection points between branches $\Sigma_1$ and
$\Sigma_2$ of the divisor.
Then $T_xX$ contains a pair of transverse symplectic two-planes, and 
the delicacy arises
since such a pair has non-trivial moduli (we have not assumed that
the planes are mutually symplectically
orthogonal).  Let $S^{\perp}$ denote the symplectic orthogonal
complement to a symplectic subspace $S$ of $T_xX$.
We can linearise the situation near $x$ and write $\Sigma_2$ as the graph of a
matrix $A: \rr^2 \rightarrow \rr^2$ over the surface
$(T_x\Sigma_1)^{\perp} \cong \rr^2$, 
viewed as a standard two-plane in $\rr^4$.    The condition that
$\Sigma_2$ be symplectic asserts that the determinant $\det(A) > -1$, whereas a
sufficient condition for a local complex model to exist (near the
transverse intersection point) is that $\det(A) > 0$.  

\vspace{0.2cm}

\noindent However, we can restore the situation with a small
perturbation near $x$ as follows.  Consider a map $u: \rr^2
\rightarrow \rr^2$ for which $\det(\Jac(u)) > -1$ everywhere and which
co-incides with the linear map $A$ outside of a large ball.  The graph
of any such $u$ defines a symplectic surface in $\rr^4$ which can be
glued back in place of $\Sigma_2$ in a neighbourhood of the point
$x$.  It follows that we want to find such a map $u$ satisfying the
additional condition that the determinant of the Jacobian of $u$ is positive at
the origin.  By choosing appropriate bases, we can suppose that the
matrix $A$ is diagonal $A = \mathrm{diag}(\lambda x, -\mu y)$ where
$\lambda > 0$ and $\mu < 0$~; the symplectic condition ensures
$\lambda \mu < 1$.  Choose some $\varepsilon > 0$ such that $\lambda \mu <
1-\varepsilon$.  We will work with deformations of the shape

$$u(x,y) \ = \ (\lambda x, f(x,y))$$

\noindent where $f(x,y) = -\mu y$ for $|(x,y)| \gg 0$ and $(\partial
  f / \partial y)|_{(0,0)} > 0$.  We can introduce a "kink" into the graph
  $(y \mapsto \mu y)$ which gives a function $f_0: \rr \rightarrow
  \rr$ with the  properties:
  (i) $f_0 (y) = 0$ only at $y=0$; (ii) $f_0(y) = -\mu y$ for $|y|$
  sufficiently large; (iii) for every $y$ the derivative
  $(df_0 / dy) > -\mu - (\varepsilon/\lambda)$; (iv) at the origin
  $(df_0 / dy)|_{y=0} > 0$.  Moreover, a suitable real one-parameter
  family of functions $(f_x)_{|x| < t}$ which interpolates between this kink
  and the linear map $(y \mapsto \mu y)$, gives
  a function $f(x,y) = f_x(y)$ 
  with $f(0,y) = f_0(y)$ and $|(\partial f / \partial y)| > -\mu -
  (\varepsilon / \lambda)$ everywhere.  For the associated map $u$, the
  graph defines a new surface $(\Sigma_2)'$ which meets $\Sigma_1$
  only at the origin $x$ of the co-ordinate system, co-incides with
  $\Sigma$ outside a ball $B$ and which at the
  intersection point is given by a linear map $\Jac (u)$ whose
  determinant is $\lambda (df_0 / dy)|_{y=0} > 0$.  Moreover, the new
  surface is again symplectic, since at every point of $B$ we see that
  $|\det(\Jac(u))| > -\lambda \mu -
  \varepsilon > -1$, where the last equality uses the fact that
  $\lambda \mu < 1-\varepsilon$ by choice of $\varepsilon$.

\vspace{0.2cm}

\noindent Hence we reduce to the situation where all the isolated
intersection points are locally given by graphs of \emph{complex}
matrices $A \in \cc^*$ (this remaining parameter is a "K{\"a}hler angle").
We can now appeal to a 
version of Weinstein's extension of the Darboux Theorem
\cite{Weinstein}.  This implies that the 
symplectic structure in a neighbourhood of $\bigcup_{i}\Sigma_{i}$ is 
uniquely determined by the areas of the $\Sigma_{i}$, their normal bundles
and the angle parameter at each intersection point.
That is, if $\Sigma'_{i}$ are symplectic surfaces in any other symplectic
 $4$-manifold $Y$,
 meeting transversally with locally positive intersections and at the
 same angles, 
such that $\Sigma'_{i}$ is diffeomorphic   to $\Sigma_{i}$, of equal 
symplectic area and with the same self-intersection numbers then 
suitable neighbourhoods of $\bigcup \Sigma'_{i}\subset Y$ and $\bigcup 
\Sigma_{i}\subset X$ are symplectomorphic. For any values of the
parameters, we may construct such a $Y$ by
``plumbing'' holomorphic line bundles over $\Sigma_{i}$, so it
suffices to prove  
the Lemma in the case when the $\Sigma_{i}$ are {\it complex curves} in a 
K{\"a}hler manifold $X$. In this case, let $\LL \rightarrow X$ be the
holomorphic line bundle 
corresponding to the divisor $\sum a_{i} \Sigma_{i}$, so there is a 
holomorphic section $\gamma$ of $\LL$ with zero-set $\sum a_{i}\Sigma_{i}$. 
The hypothesis asserts 
that the degree of $\LL$ over each Riemann surface $\Sigma_{i}$ is
nonnegative.  
It is clear then that we can find a  $C^{\infty}$ section $s$ of $\LL$ over a 
neighbourhood of $\bigcup_{i} \Sigma_{i}$ with the following properties.
                       
\begin{enumerate}

\item The section $s$ is holomorphic and non-vanishing in a neighbourhood of 
the intersections $\Sigma_{i}\cap \Sigma_{j}$. In fact we can arrange that $s$ 
defines a local holomorphic trivialisation of $\LL$ such that, in suitable 
local complex co-ordinates $z_{1},z_{2}$, the section $\gamma$ is
$\gamma = z_{1}^{a_{i}} z_{2}^{a_{j}} s.$

\item The section $s$ is transverse to zero and
 holomorphic in a neighbourhood of its zero set. In fact we can arrange that
around any zero of the section $s$ on  $\Sigma_{i}$ there is a local 
holomorphic trivialisation of $\LL$ and local complex co-ordinates in which
$s$ and $\gamma$ are  represented by the functions
$z_{2}$ and $z_{1}^{a_{i}} $ respectively.

\end{enumerate}

\noindent Now, for small non-zero $\epsilon$, let $\Sigma_{\epsilon}$
be the zero-set of the 
section $\gamma- \epsilon s$ of $\LL$. It is clear that, if $\epsilon$
is small enough, $\Sigma_{\epsilon}$ is 
compact and, provided it is cut out transversally, represents the homology 
class $\sum a_{i}[\Sigma_{i}]$. We need to prove that $\Sigma_{\epsilon}$ is 
a symplectic surface, cut out transversally. To do this it suffices to show 
that for any point $p$ of $\bigcup_{i} \Sigma_{i}$ there is a neighbourhood
$N_{p}$ of $p$ and a $\delta_{p}>0$ such that $\Sigma_{\epsilon}\cap N_{p}$ 
is a symplectic surface, cut out transversally, provided 
$0< | \epsilon | <\delta_{p}$. 
Relabelling the $\Sigma_i$ appropriately, there are three cases to consider:

\begin{itemize}

\item $p$ is an intersection point of  two surfaces, $\Sigma_{1}$ and
  $\Sigma_{2}$;
\item $p$ lies on a single surface $\Sigma_{1}$, and $s$ does not vanish 
at $p$;
\item $p$ is a zero of $s$ on $\Sigma_{1}$.

\end{itemize}

\noindent The first and third cases are immediate. In these cases the
section $s$ is  
holomorphic near $p$ and $\Sigma_{\epsilon}$ is locally a complex variety 
defined by  equations

$$ z_{1}^{a_{1}} z_{2}^{a_{2}} = \epsilon  \qquad  \mathrm{or} \qquad 
    z_{1}^{a_{1}} = \epsilon  z_{2}, $$
                                    
\noindent respectively.
Since these equations cut out transversally smooth curves in $\cc^{2}$, 
for any non-zero $\epsilon$, the result follows.
In the second case, we can choose local complex co-ordinates $z_{1},z_{2}$ and 
a local holomorphic trivialisation of $\LL$ so that $\gamma$ is represented by 
the function $z_{1}^{a_{1}} $, while $s $ is represented by 
a smooth function $F(z_{1},z_{2})$,
with $F(0,0)= C\neq 0$.  If we put $\eta= \epsilon^{1/a_{1}}$,
then $\Sigma_{\epsilon}$ is locally defined by an equation
$$ z_{1}^{a_{1}} = \eta^{a_{1}} F(z_{1},z_{2}).$$ 
Now let $w_{1},w_{2}$ be rescaled co-ordinates $w_{i}= \eta^{-1} z_{i}$. The 
equation becomes
$$     w_{1}^{a_{1}}  = F(\eta w_{1},\eta w_{2}). $$
We view this as a perturbation of the equation $\{ w_{1}^{a_{1}} = C \}$ whose 
solution set is the union of $a_{1}$ complex lines parallel to the $w_{2}$ 
axis. It is then easy to see that the zero-set of the perturbed equation has 
the desired properties over a ball  $\Vert (w_{1},w_{2})\Vert \leq 
r\eta^{-1}$ provided that $r$ 
and $\eta$ are sufficiently small, so we can take $N_{p}$ to  be the ball
$\Vert (z_{1},z_{2})\Vert < r$. 
\end{Pf}

\noindent Note that the same proof can be extended to obtain a more
general result, in which the
surfaces $\Sigma_{i}$ are allowed to have singularities modelled on singular 
points of complex curves.  A particular variant we shall use allows us to
separate out $(-1)$-sphere components: here is the simplest case.

\begin{Lem} \label{smoothingtwo}
Let $C+2E$ be a positive symplectic divisor comprising an embedded
symplectic sphere $E$ of square $(-1)$ with multiplicity $2$ and
a reduced component $C$ which meets $E$ transversely once.  Then in an
arbitrarily small neighbourhood of $C \cup E$ there is a smooth
symplectic surface which comprises a
multiplicity one copy of $E$ and a disjoint surface
representing $C+E$.
\end{Lem}

\begin{Pf} By plumbing, we can again assume that the surfaces $C$ and
  $E$ are complex curves in a K{\"a}hler surface.  The divisor $C+2E$
  is cut out by a section $\nu$ of a line bundle $\mathcal{L}'$;
  if we take the reduced divisor $(C+2E)_{red}$ we obtain a
  holomorphic section $\eta$ of the line bundle $\mathcal{L}$ with
  first Chern class $C+E$.  Observe that $(C+E) \cdot E = 0$.  We can
  choose a smooth section $s$ of
  $\mathcal{L}$ which is nowhere vanishing on $\mathcal{L}|_E$, which
  is holomorphic near its zeroes and which near the intersection point
  $p$ of $C$ and $E$ trivialises $\mathcal{L}$ and such that in local
  co-ordinates the section $\eta$ is given by $\eta = z_1 z_2^2 s$.
  Then, for sufficiently small $\epsilon$, the zero-set of $\eta -
  \epsilon s$ will be a symplectic surface of the required structure
  by arguments as above.
\end{Pf}

\noindent Again there is an obvious generalisation, allowing one to
  separate out a number of embedded exceptional curves which may
  themselves have higher multiplicities.

\vspace{0.2cm}

\noindent We close the section with a (perhaps helpful) analogy.  It
follows from (\ref{lefpencildefn}) that we can always
choose a compatible almost-complex structure $j$ on $X'$ such that
\begin{itemize}
\item the projection map to $\sss^2$ is pseudoholomorphic;
\item the structure is integrable in a sufficiently small tubular
  neighbourhood of each singular fibre, i.e. over discs
  $D_{\varepsilon} \subset \sss^2$ centred on the $f(p_i)$.
\end{itemize}
\noindent We thus induce a map
from the smooth part of the base $\sss^2 \backslash \{ f(p_1), \ldots,
f(p_r) \}$ to the Deligne-Mumford coarse moduli space of curves $M_g$;
from the local models, this
map extends to a map of the closed sphere $\phi_f: \sss^2 \rightarrow \mgbar$
into the projective moduli space of \emph{stable} curves.  (For a
generic choice of $j$, and a pencil of curves of genus $g>3$, the
sphere $\phi_f (\sss^2)$ will be disjoint
from the orbifold loci of curves with automorphisms.  Away from these
loci the moduli space is fine, and if we wish we can
invoke the existence of a ``universal curve''.)  This sphere
has transverse locally positive intersections with the divisor of
stable curves;  transversality provides nodal singular fibres, and the
positivity shows the local complex co-ordinates have the correct
orientations to construct, directly from the topology, a symplectic
structure on $X'$.
Similarly, standard surfaces can be obtained from
sections of symmetric product fibrations which have transverse locally
positive intersections with natural diagonal strata: the
transversality ensures the surfaces are smooth, and the positivity fixes
the orientations to be compatible with the symplectic structure
produced above.


\section{The Abel-Jacobi map}

The rest of the proof shall involve the relationship between three
fibre bundles which can be associated to a smooth surface fibration.
We begin by reviewing some of the standard theory of  ``divisors''
on  a  Riemann surface $S$. We write $\Sym^{r}(Z)$ for the $r$-fold
``symmetric product'' of a space $Z$: the quotient of the product of $r$ 
copies of $Z$ by the action of the permutation group on $r$
letters. Elements of $\Sym^{r}(Z)$
can also be regarded as formal sums $ \sum a_{i} p_{i}$ 
where the $a_{i}$ are positive integers with $\sum  a_{i} = r$ and the
$p_i$ are points of $Z$. The
starting point  
is the fact that $\Sym^{r}(\cc)$ can be identified with
$\cc^{r}$ via the elementary symmetric functions. An $r$-tuple
$(z_{1},\dots, z_{r})$ maps to  $\usigma(z_{1},\dots z_{r})= 
(\sigma_{1},\dots \sigma_{r})$ where 

\begin{Eqn} \label{induced}
\sigma_j = \sum_{u \in \Sym(j)} \, \prod_{i=1}^j  z_{u(i)}; \qquad \mathrm{so} \qquad \sigma_{1} = \sum z_{i},
  \ \  \sigma_{2} = \sum_{i \neq j} z_{i} z_{j}, \ \ldots, \ \sigma_r =
  z_1 \cdots z_r.  
\end{Eqn} 

\noindent Suppose we have a holomorphic diffeomorphism
$\alpha:U\rightarrow V$ between  
open sets $U,V$ in $\cc$.  Write, for any such $U$, $\usigma(U)$ for
the image of $\usigma: \Sym^r(U) \rightarrow \cc^r$.  Then we get an
induced holomorphic  diffeomorphism
$$\Sym^{r}(\alpha):\usigma(U)\rightarrow \usigma(V). $$
It follows from this that if $S$ is any Riemann surface
the symmetric product $\Sym^{r}(S)$ has a natural structure of an
$r$-dimensional complex manifold.  Now suppose $S$ is compact of genus $g$.
The Picard variety
$\Pic_{r}(S)$  of line bundles of degree $r$ is a $g$-dimensional complex torus
which, up to a choice of origin, can be identified with 
$$\Pic_{0}(S) = H^{1}(S;\mathcal{O})/ H^{1}(S;\zz). $$
The complex vector space $H^{1}(S;\mathcal{O})$ is the dual of the space of 
holomorphic $1$-forms $H^{0}(\Omega^{1}_{S}) = H^0 (K_S)$.
There is a canonical holomorphic map - the \emph{Abel-Jacobi map} - 

$$\tau_{r}: \Sym^{r}(S) \ \rightarrow \ \Pic_{r}(S), $$

\noindent which can be thought of either as the map assigning a line
bundle to a divisor  
or as the map induced by integration of  holomorphic $1$-forms along paths. 
The fibres of $\tau_{r}$ are complex projective spaces: the fibre over
$\mathcal{L}$ is $\pp (H^{0}(S; \mathcal{L}))$.  For future reference,
here is an easy but important topological lemma, which will play a
role in our considerations of bubbling (\ref{bubblingfacts}): a proof
can be found in \cite{Bertram-Thaddeus}. 

\begin{Lem} \label{secondhomotopy}
For every $r>1$ (and every genus) the Hurewicz homomorphism  
$$\pi_2 (\Sym^r (\Sigma)) \rightarrow H_2 (\Sym^r (\Sigma))$$ 
\noindent has rank precisely one.
\end{Lem}

\noindent Now consider in particular the case when $r=2g-2$.  The
canonical line bundle 
$K_{S}$ gives a preferred point\footnote{We will sometimes refer to this
  point as ``zero'', and the section of a bundle of Picard varieties
  defined by the points representing the canonical line bundles as the
  ``zero-section''.} in
$\Pic_{2g-2}$.  The fibre of $\tau_r$
over $K_{S}$  
is a copy of $\cc \pp^{g-1}$, whereas all the other fibres are copies of 
$\cc \pp^{g-2}$. Thus, away from $K_{S}$, the map $\tau_{2g-2}$ is the
projection  
of a holomorphic $\cc \pp^{g-2}$ bundle;  in general, at a divisor
$p_1 + \cdots + p_r$ the image of $d\tau_r$ is dual to the subspace 

$$ (\im (d\tau_r)_{\sum p_i})^* \ = \ \{ \nu \in H^0(\Sigma,
K_{\Sigma}) \ | \ \mathrm{div}(\nu) \geq p_1 + \cdots + p_r \}.$$

\noindent Here $\mathrm{div}(\nu) \geq D$, by definition of notation,
if the zero-set of $\nu$ contains the divisor $D$.  Hence, if $\sum
p_i$ is not canonical then necessarily $\nu \equiv 0$ and $d\tau_r$ is
surjective: if $\sum p_i = \mathrm{Zeroes}(\nu)$ for $\nu \in H^0 (K_{\Sigma})$
then $d\tau_r$ has a one-dimensional cokernel. 

\vspace{0.2cm}

\noindent Families of elliptic operators can be understood through
their linearisations, and the structure of the map around the 
special fibre can be identified with a standard model. Let $V$ be a complex
vector space of dimension $g$ and let $M$ be the subet of $V^{*}\times \pp(V)$
defined by

\begin{Eqn} \label{localmodel} 
M= \{ (\theta,[x]): \theta(x) = 0 \}.
\end{Eqn}

\noindent Let $A:M\rightarrow V^*$ be the projection map onto the first
factor. The  
fibre of $A$ over $0$ is $\pp(V)\cong \cc \pp^{g-1}$,
while the fibre over a non-zero $\theta$ is $\pp({\ker}(\theta))
\cong \cc \pp^{g-2}$. Now take $V=H^{0}(S; K_S)$ and identify a
neighbourhood of $K_{S}$ 
in $\Pic_{2g-2}(S)$ with a neighbourhood of $0$ in $V^{*}$ via the
canonical flat  
structure on the torus.

\begin{Lem}  \label{localmodelforAJmap}  Let $S$ be a compact Riemann
  surface;  write $r=2g-2$ as usual.  Then
$(\Sym^{r}(S), \tau_{r})$ is holomorphically equivalent to $(M,A)$ in
neighbourhoods of the fibres over $K_S$ and over $0$ respectively. 
\end{Lem}
 
\begin{Pf}
Points in a neighbourhood of $0 \in V^*$ correspond to line bundles
whose holonomies around loops are all small; we can therefore regard
the holonomy of such a point as a real-valued, rather than
circle-valued, one-form $a$.  (Via
``holonomy difference'', we can parametrise a neighbourhood of $K \in
\Pic_r(S)$ in the same fashion.)   On $S$ there is a canonical
isomorphism between the real vector space $H^1 (S; \rr)$ and the
complex vector space
of $(0,1)$-forms $H^{0,1}(S)$.  Under this isomorphism, $a$ maps to a
form $\nu$ which defines a local holomorphic parameter on the Jacobian
and hence Picard varieties.  To see this directly, regard the space of
holomorphic 
bundles as a space of unitary connexions on a fixed Hermitian line
bundle.  Two such connexions differ precisely by a $(0,1)$-form,
unique up to gauge; in a small disc around $\cdbar_K$, the d-bar
operator defined by the canonical line bundle, we can fix the gauge
consistently.  It follows that $\cdbar_K +\nu$ is the 
$\cdbar$-operator on the holomorphic line bundle defined by $a \in V^*$.
The holomorphic
sections $\phi$ of $\mathcal{L}_a$  by definition satisfy 

$$\cdbar_a (\phi) = \cdbar_K \phi + \nu \wedge \phi = 0.$$

\noindent We have an embedding into this space of those sections $s
\in H^0 (S;K)$ which satisfy (by 
definition) $\cdbar_K (s) = 0$ and for which $\langle \nu, s \rangle =
0$. Since the two vector spaces have the same dimension, this
embedding is an isomorphism, and identifies the linear system of
sections of $\mathcal{L}_a$ with the fibre of the map $A$.  Since the
(scheme-theoretic) fibres of $\tau_r$ are
precisely linear systems, the result follows.
\end{Pf}

\noindent Return to the Lefschetz fibration $X' \rightarrow \sss^2$.
Fix an almost-complex structure $J$ on $X'$ for which the projection
map is holomorphic; write $X^*$ for the complement in $X'$ of the
singular fibres. The 
restriction of $f$ to $X^{*}$ is a genuine $C^{\infty}$ fibration and each 
fibre has the structure of  a Riemann surface. We define the fibrewise 
symmetric product $\Sym^{r}_{f}(X^{*}) = X^*_r (f)$, 
as a topological manifold, in the
obvious way; there is a projection map from $X^*_r (f)$ to 
$\sss^{2}$ whose fibres are the symmetric products of the fibres of 
$f$. (Since a homeomorphism of a space defines a homeomorphism of its 
symmetric products the construction at this level can just be seen as the 
ordinary construction of an associated bundle, associated in this case to the 
group of self-homeomorphisms of a surface.)

\vspace{0.2cm}

\noindent To put a smooth manifold structure on $X^*_r(f)$
requires more care,  
since a diffeomorphism of a surface does \emph{not} induce diffeomorphisms of 
its symmetric products. 

\begin{Defn} \label{restrictedcharts}
A \emph{restricted chart} on $X^{*}$ comprises 
a diffeomorphism $\chi:D\times D\rightarrow X^{*}$ (where $D$ is the unit disc 
in $\cc$) such that
\begin{enumerate}
\item there is a holomorphic diffeomorphism $\theta$  from $D$ to an
  open set in $\sss^{2}$ such that the diagram below commutes:

\[\xymatrix{
D \times D \ar[r]^{\chi} \ar[d]_{\mathrm{pr}_1} & X^* \ar[d]^f \\
D \ar[r]^{\theta} & \sss^2
}\]

\item for each fixed $\tau\in D$ the map $\chi_{\tau}$ defined by 
$\chi_{\tau}(z) = \chi(\tau,z)$ gives a holomorphic diffeomorphism from $D$ to 
an open set in the fibre $f^{-1}( \theta(\tau))$, where the latter is
endowed with complex structure coming from $J$ on $X'$. 
\end{enumerate}
\end{Defn}

\noindent It is true, but not completely trivial, that any point of
$X^{*}$ lies in the  
image of  a restricted chart. (This follows from the Riemann mapping theorem, 
as in \cite{AhlforsBersRMT} for example, with smooth dependence on
parameters.) Two restricted  
charts compare, on the intersection of their images, by a smooth family of 
holomorphic diffeomorphisms of open sets in $\cc$. Since these
holomorphic  
diffeomorphisms induce holomorphic diffeomorphisms of the symmetric products 
as after (\ref{induced}) above, and the introduction of smooth parameters
introduces no  
difficulties, these restricted charts induce a system of charts on 
$X^*_r(f)$ which differ on their overlaps by smooth maps. In sum we 
define $X^*_r(f)$ as a smooth manifold, with a smooth fibration over
$\sss^{2}\backslash\{ f(p_{i})\}$. Each fibre in $X^*_r(f)$
has a natural  
complex structure, and these vary smoothly in the obvious sense. 

\vspace{0.2cm}

\noindent All of this discussion is compatible with the construction
of the Picard  
varieties. We now specialise to the case when $d=2g-2$, where $g$ is the genus 
of the fibres of $f$.  There is one subtlety, which we introduce now
and clarify at the end of the section.  If $Z \rightarrow B$ is
any smooth fibration of curves, with smoothly varying holomorphic
structures on the fibres in the above sense, there is a vector bundle
$W \rightarrow B$ whose fibres are \emph{canonically} identified with
the spaces of holomorphic one-forms on the fibres.  On the other hand,
if $Z$ is holomorphic then $V=f_* K_Z$ and $W$ are \emph{not equal}, but
differ by twisting by the canonical bundle of the base:  $f_* K_Z = W
\otimes K_B$.  In constructing submanifolds representing the homology
class $K_Z$, for a given manifold $Z$, it is sections of this twisted
bundle $V$ that will be of primary importance. 
Thus, for the smooth part of a Lefschetz fibration $f: X' \rightarrow
\sss^2 \backslash \{ \Crit \}$, we have:

\begin{enumerate}
\item A rank $g$ complex vector bundle $W$ over
$\sss^{2} \backslash \{ f(p_{i})\}$ whose fibres are canonically
identified with 
the holomorphic  $1$-forms on the fibres of $f$.
\item A bundle of complex tori $\Pic_r^f(X^{*})=P^*_r(f)$
 over $\sss^{2}\setminus\{f(p_{i})\}$,
with a ``zero-section'', 
whose fibres are quotients of the fibres of $W^{*}$ by the integer lattices
$H^{1}(f^{-1}(\tau);\zz)$. 
\item A smooth map $\tau_{f}: X^*_r(f) \rightarrow P^*_r(f)$,
  commuting with the projection 
maps, which is equal to the Abel-Jacobi map $\tau_{r}$ on each fibre.
\item A diffeomorphism, compatible with the projection maps,
 between  a neighbourhood of $(\tau_{f})^{-1}(0)\subset 
X^*_r(f)$ and a neighbourhood of $A^{-1}(0)$ in
$\underline{M}(V)$,  
which agrees with the map considered above on each fibre.
\end{enumerate}

\noindent Here $\underline{M}(V)$ is the space, fibring over $\sss^2
\backslash \{ f(p_i) \} =\Delta$, constructed from the
vector bundle $V = W \otimes K_{\Delta}$ in  
the obvious way.  Precisely, inside the fibre product of $W^*
\rightarrow \Delta$ and $\pp(V) \rightarrow \Delta$ we take the subvariety
defined by (\ref{localmodel}) in each fibre\footnote{This is
  well-defined, and smoothly equal to $V^* \times_{\pi} \pp(V)$, as
  the spaces $V$ and $W$ are isomorphic up to scale.}.  Note that a
$\cdbar$-operator (for instance from a connexion) on
the complex vector bundle $V$ will induce an integrable holomorphic
structure on the total space of this model $\underline{M}(V)$.

\vspace{0.2cm}

\noindent Before turning to the singular fibres, we need one more
ingredient.  Recall that  we have a collection of exceptional curves
$E_{i}$ in $X'$. Let $\Lambda \rightarrow X'$ be the 
corresponding line bundle, so $c_{1}(\Lambda)$ is the Poincare dual of the 
sum of the exceptional curves $E_{i}\subset X'$. This line bundle 
induces a line bundle $\Lambda^*_{r}(f)$ over $X^*_r(f)$. To
obtain this as  
a topological line bundle we just take the obvious lift of the action of the 
symmetric group on the tensor product
$$  \pi_{1}^{*}(\Lambda)\otimes \cdots \otimes \pi_{r}^{*}(\Lambda)$$ 
over the fibrewise product of $r$ copies of $X^{*}$ (where $\pi_{i}$ are
the projection maps to each factor).  One can check, although we do not need 
this, that this line bundle has a natural smooth structure.  

\vspace{0.2cm}

\noindent We now extend each of the ingredients above to the nodal
curves in the Lefschetz fibration.
Recall that we have chosen 
a genuine complex structure on $X'$ in the neighbourhood of each singular 
fibre, and the map $f$ is holomorphic. Thus the discussion can take place 
entirely in the realm of complex geometry.  We need, then, a suitable
extension of the discussion of divisors and line bundles to singular 
curves.  Such a theory is well established in algebraic geometry;  we
will give a ``users' guide'' to this in the Appendix, quoting the
substantial theory from the literature.   For now we simply summarise  the 
result that we need.

\begin{Thm} \label{itallextends}
Suppose $f: \mathcal{X} \rightarrow D$ is a holomorphic map from a
smooth complex  
surface to the disc such that the fibres over $t \neq 0$ are smooth
but  $f^{-1}(0)$ is an irreducible curve with a single quadratic
singularity $Q$.  
Let $E_{i}$ be sections of $f$. Then there are holomorphic vector 
bundles $W$ and $V = W \otimes K_D$ over $D$, and 
smooth complex manifolds $\mathcal{X}_r(f)$ and $\mathcal{P}_r(f)$,  with
commuting maps 

\[\xymatrix{
\pp(V) \ar[r]^{\iota} \ar[d] & \mathcal{X}_r(f) \ar[r]^{\tau}
\ar[d]_F & \mathcal{P}_r(f) \ar[d]^{\pi} \\
D \ar[r]^{=} & D \ar[r]^{=} & D
}\]

\noindent and a relatively ample
holomorphic line bundle $\Lambda_r (\pi)$ over $\mathcal{X}_r(\pi)$ 
such that:

\begin{itemize}
\item when restricted to the punctured disc the data agrees with the 
symmetric poduct construction above, giving on each fibre the maps

$$\pp (H^0 (K_{\Sigma})) \stackrel{\iota}{\longrightarrow} \Sym^r(\Sigma)
\rightarrow \Pic_r (\Sigma);$$

\item the line bundle $\Lambda_r (\pi)$ restricts on each smooth fibre
  to the line bundle defined by the sum of the divisors in $\Sym^r
  (\Sigma)$ of $r$-tuples of points meeting one of the $E_i$;

\item the spaces $(\mathcal{X}_r(f))_0$ and $(\mathcal{P}_r(f))_0$
  over $0$  are irreducible and 
have simple normal crossing singularities:
in suitable co-ordinates around the singularities 
the vertical maps in the diagram above are given by the standard model 
$(z_{1},\dots z_{n})\rightarrow z_{1} z_{2};$

\item the image of the map $\iota$  does not meet the set of critical
  values of the map $F$ and 
there is a holomorphic diffeomorphism
between a neighbourhood of this image and a neighbourhood of the exceptional 
fibres in the complex manifold $\underline{M}(V)$;

\item the canonical isomorphism, over $f^{-1}(D\setminus\{0\})$,
between the tangent space along the fibres 
of $\pi$ and the pull-back of the bundle $W^{*}$,
extends to the smooth part of $f^{-1}(0)$; 

\item a Zariski-open subset of $F^{-1}(0)$ can be identified with the 
open complex manifold $\Sym^{r}(\mathcal{X}_{0}\setminus \{Q\})$.   
\end{itemize}
\end{Thm}

\noindent Using this we can immediately extend our constructions over
the whole space  
$X'$. We obtain a vector bundle $V\rightarrow \sss^{2}$, spaces $X_r
(f)$, $P_r (f)$  and  a line 
bundle $\Lambda \rightarrow X_r (f)$, with a diagram of maps 

$$   \Lambda\rightarrow X_r(f)  \rightarrow P_r (f) \leftarrow
\underline{M}(V), $$

\noindent where all the data fibres over $\sss^{2}$.  Although the
difference between the bundles $V$ and $W$ seems rather artificial
above, where the canonical bundles of the (open) base surfaces were
trivial, it now plays an important role.  There is an embedding
$\pp(V) \rightarrow X_r(f)$ but a neighbourhood of the section of
$P_r (f)$ defined by the canonical lines is naturally identified
with a neighbourhood of zero in the bundle 

$$W^* = (V \otimes K_{\sss^2}^{-1})^* = V^* \otimes \mathcal{O}(-2).$$

\noindent Thus although $\pp(V) = \pp(W)$,  sections of the
projective bundle arising from nowhere zero sections of $V$ and $W$
lie in different homotopy classes.


\section{Symplectic structures and the Gromov invariant}

In order to count pseudoholomorphic sections of the map $X_r(f)
\rightarrow \sss^2$ we will need to discuss almost-complex and
symplectic structures on the total space.  
We begin with a general discussion of almost-complex structures on fibre 
bundles. Suppose we have an exact sequence of $\rr$-linear maps

\begin{Eqn} \label{sequence}
0\rightarrow  U\stackrel{\iota}{\longrightarrow} V \longrightarrow W
\rightarrow 0 
\end{Eqn}

\noindent where $U$ and $W$ are complex vector spaces but $V$ is
initially only a real  
vector space. We want to consider the set $\mathcal{S}$ of complex
structures on 
$V$ such that the  
given maps are complex linear. We can regard a complex structure on $V$ as
a complex subspace $V'$ of the complexification $V_{\cc} = V\otimes_{\rr} \cc$ 
such that $V_{\cc} = V'\oplus V''$, where $V''$ is the conjugate of $V'$.
From this one sees that the compatible complex structures correspond to {\it 
splittings} of the sequence of complex linear maps

$$U \ \rightarrow \ V_{\cc}/ \iota(\overline{U}) \ \rightarrow \
\overline{W}, $$ 

\noindent obtained from (\ref{sequence})  by complexifying and then
writing $U_{\cc}= U\oplus \overline{U}$ and $W_{\cc}=
W\oplus\overline{W}$. This means that there is a natural {\it affine} 
structure on $\mathcal{S}$: two elements of $\mathcal{S}$ differ by a map in
$\Hom( \overline{W}, U).$

\vspace{0.2cm}

\noindent Now consider a manifold $Z$  which is the total space of a
smooth fibre bundle 
$g:Z\rightarrow B$ and suppose that we  are given almost-complex
structures on the  
base $B$ and on each fibre of $g$, 
varying smoothly in the obvious sense. Then the compatible almost-complex 
structures on $Z$ are sections of a bundle over $Z$ whose fibres are affine 
copies of $\Hom(g^{*}TB, T^{vt}Z)$, where $T^{vt}Z$ denotes the
``vertical'' tangent spaces 
along the fibres of $g$.  We see at once from this that such almost-complex 
structures {\it exist}, since we can take affine linear combinations using a 
partition of unity. More generally if we are given a structure over an open 
set $U\subset Z$ and we have a closed set $K\subset U$ we may find a structure 
over the whole of $Z$ which agrees with the given one over a neighbourhood of 
$K$.  Taking $K$ to be some closed neighbourhood of the singular
fibres of $X_r(f) \rightarrow \sss^2$, equipped with its integrable
complex structure, gives the

\begin{Lem} \label{almostcomplex}
Let $\mathcal{J}$ denote the class of smooth almost complex structures 
on $X_r(f)$ which agree with the
given integrable structures on the fibres and in some (not necessarily
fixed) open neighbourhoods of the 
singular fibres, and for which the projection map to the sphere is
holomorphic.  Then $\mathcal{J}$ is a non-empty smooth
infinite-dimensional manifold.  
\end{Lem}

\noindent  We should stress that we will always deal with almost
complex structures drawn from this class.  To make this discussion
more concrete, return to $g: Z \rightarrow B$ and
suppose  that the base $B$ has complex dimension $1$; in fact 
let us suppose that $B$ is an open set in $\cc$. Then two compatible 
almost-complex structures differ by a smooth vector field along the fibres.
Explicitly, suppose we have chosen preferred co-ordinates 
$(z_{i},\tau)$ in a 
neighbourhood of a point of $Z$ in the manner of
(\ref{restrictedcharts}). So the $z_{i}$ are 
holomorphic co-ordinates along the fibres and $\tau$ is the canonical 
co-ordinate in the base.  In these co-ordinates we have a distinguished 
almost-complex structure given by the identification with an open set
in $\cc^n$. Any other almost-complex structure is represented by a vector field
$v_{i}(z_{j};\tau)$. For example the $\overline{\partial}$-operator of this  
structure is just

$$   \overline{\partial}_{v}(f) = \left(
 \frac{\partial f}{\partial \overline z_{i}}; \ \frac{\partial f}{\partial 
\overline{\tau}} + \sum v_{i} \frac{\partial f}{\partial \overline
z_{i}} \right). $$

\noindent Similarly a section $z_{i}= \phi_{i}(\tau)$ defines a
pseudoholomorphic curve in this structure if
 
$$\frac{\partial \phi_{i}}{\partial \overline{\tau}} = -\overline{v}_{i}. $$

\noindent If our structure is given by a vector field $v_{i}$ in one
set of co-ordinates, and we consider a set differing by a family of 
holomorphic maps $g_{\tau}$, the same structure is represented in the
new co-ordinates by 
   
$$(Dg_{\tau}) v + \frac{\partial g_{\tau}}{\partial \overline{\tau}}.$$

\noindent We will use these explicit presentations later to obtain
almost-complex structures on $X_r(f)$ with helpful ``genericity''
properties.  For now, we turn to symplectic structures. According to
Gompf \cite{GompfGokova}, given any map of compact almost-complex
manifolds $\pi: (Z; J_Z) \rightarrow (B;J_B)$  for which

\begin{itemize}
\item $B$ admits a symplectic structure;
\item for each $b \in B$ there is a neighbourhood $W_b$ of
  $\pi^{-1}(b)$ in $Z$ with a closed two-form $\eta_b$ which tames the
  almost-complex structure $J_Z |_{\ker d\pi}$;
\item the $\eta_b$ are all induced by a single non-zero cohomology
  class $c \in H^2 (Z; \zz)$,
\end{itemize}

\noindent then $Z$ admits a symplectic structure.  The proof uses a
patching argument to obtain vertical non-degeneracy and Thurston's
trick of pulling back a large multiple of the base form to obtain
global non-degeneracy.  In our situation, the subvariety of
the symmetric product of a Riemann surface consisting of divisors
containing a given point 
$e$ is dual to a K{\"a}hler form.  (To see that line bundle is ample,
one can use the Nakai-Moishezon criterion, for instance.)  This
determines a class $c$ and the
positivity properties of the bundle $\Lambda_r(\pi)$ in
(\ref{itallextends}) take care 
of the singular fibres, providing a K{\"a}hler form near the ends with
which we can patch the symplectic forms on $X^*_r(f)$.  Putting this
together, we obtain symplectic
forms $\Omega$ which have the shape 

\begin{Eqn} \label{symplonsymm}
\Omega \ = \ \mu_{\Sym} + R f^* \omega_{st}; \qquad R \gg 0
\end{Eqn}

\noindent where $\mu_{\Sym}$ is some closed two-form on $X_r(f)$
which is symplectic on the fibres and $\omega_{st}$ is the standard
symplectic form on $\sss^2$ as before.  Note that all such symplectic
forms $\Omega$ are deformation equivalent, once we have chosen $c$
(which for us comes from the exceptional sections of $f$).  It is easy
to see that a given
almost-complex structure $J \in \mathcal{J}$, as  provided by
(\ref{almostcomplex}), is tamed
by the forms (\ref{symplonsymm}) once $R > R(J)$ is sufficiently large
(where the precise value will depend on $J$).

\begin{Lem} \label{easyforms}
$X_r (f)$ admits
symplectic structures $\Omega$ which restrict on each fibre to the usual
K{\"a}hler structure induced from the integrable complex structures on the
fibres of $f$.   
\end{Lem}

\noindent It seems likely that a more precise statement holds.  For a
space $Z$ equipped
with a map $Z \rightarrow \sss^2$ let $\Gamma(Z)$ denote the
set of homotopy classes of sections.  We have natural maps

\begin{Eqn} \label{naturalmaps}
H^2 (X'; \zz) \ \stackrel{\mu}{\longrightarrow} \ H^2
(X_r(f);\zz);
\end{Eqn}
\begin{Eqn} \label{naturalmaps2}
\Gamma (X_r(f)) \ \stackrel{\nu}{\longrightarrow} \ H_2 (X'; \zz)
\end{Eqn}

\noindent defined as follows.  For the first, represent the class $c$
as the first Chern class of a line bundle $L_c \rightarrow X'$.  This
defines a topological line bundle on $X_r (f)$ whose fibre at a tuple
$p_1 + \cdots + 
p_r$ is precisely the quotient of $\otimes_i (L_c)_{p_i}$ by the
symmetric group; take the first Chern
class to define $\mu$.  For the second, choose a smooth
section $\psi$ in the homotopy class and define $\nu(\psi)$ by taking
the $r$-tuples 
of points in each 
fibre $f^{-1}(t)$ of $X'$ designated by the value $\psi(t)$.  (Thus
$\psi$ defines a closed subset $C_{\psi}$ of $X'$ and $\nu(\psi) =
[C_{\psi}]$).  Then we have the

\begin{Question} \label{naturalform}
If $\omega$ is an integral symplectic form on $X'$ does $\mu(\omega)$
contain symplectic forms on $X_r(f)$? 
\end{Question}

\noindent If this is true, one would expect that for a section $\psi
\in \Gamma(X_r(f))$, we have an identity
$\langle \mu(\omega), [\psi] \rangle \ = \ \langle \omega,
[\nu(\psi)] \rangle$.  This would give another approach to the
question of \emph{bubbling} that we will tackle, by different means,
later on.  (The reason is simply that identifying a particular
symplectic form on $X_r(f)$ enables us to
estimate the area of bubbles in terms of data on $X'$.)  Before
computing the Gromov invariant of relevance to us,
we recall the general theory. There is now a rigorous definition of
Gromov-Witten invariants, valid for any taming almost-complex
structure $J$ on any closed symplectic manifold and independent of
regularity or monotonicity hypotheses.   (Invariants counting
holomorphic sections were also considered by Seidel \cite{Seidel2}, who
arranged the invariants for different homology classes of section into
an element of a suitable Novikov ring.)  One can proceed in two directions:
develop a theory for ``generic'' almost-complex structures, and prove
a lemma on how to compute the invariants given a non-generic structure
for which the moduli space is compact and smooth but of the wrong
dimension (cf. \ref{regularspaces}), or develop a theory via
``virtual classes'' which works at once for any almost-complex
structure, at the cost of laying heavier analytic foundations.
Although our computations use obstruction bundles, we are in the former
``elementary'' case above, and our proof does not require analysis
beyond that presented in \cite{Gromov} to produce the desired
invariant. 

\vspace{0.2cm}

\noindent Recall $\mathcal{J}$ denotes the set of
smooth almost-complex structures on $X_r(f)$ described in
(\ref{almostcomplex}).  Write $[\psi_V]$ for the
homology class of sections of $X_r(f)$ defined by a smooth nowhere
zero section of the
vector bundle $V \rightarrow \sss^2$ and the embedding $\pp(V)
\rightarrow X_r(f)$.  (That such sections exist is trivial if the
degree of the pencil $k$ is large, and hence the rank $g$ of $V$ is at
least three (\ref{highdegree}).)  To define the Gromov invariant, we
will need to understand bubbling, and so we begin with the relevant
technical Lemma:  suppose our Lefschetz pencil is by surfaces of genus $g>2$.

\begin{Lem} \label{bubblingfacts}
Two-spheres inside the fibres of $F: X_r(f) \rightarrow \sss^2$ are
governed by the following constraints:

\begin{enumerate}
\item The second homotopy group $\pi_2 (\Sym^{2g-2}(\Sigma)) = \zz$.  A
generator $h$ for this group can be given by a line - a rational curve
of degree 1 - inside a projective space fibre of the Abel-Jacobi map. 

\item Let $(\phi_j)$ be a family of holomorphic sections of $F$ in the
homology class $[\psi_V]$, and suppose that the $\phi_j
\rightarrow \phi_{\infty}$ converge to a cusp curve in the sense of
Gromov.  Each bubble component is homologous
to a multiple of $h$.

\item To each bubble $\phi$ in a fibre of $X_r(f)$ we can
  associate a closed subset
  $C_{\phi}$ of $X'$ (cf. \ref{naturalmaps2}).  If the bubble lies in
  the homology class $Nh$ then $[C_{\phi}] = N[\Fibre]$.
\end{enumerate}

\end{Lem}

\begin{Pf}
The first statement is a sharper version of (\ref{secondhomotopy}),
and can be checked for instance by  mapping $\Sym^{2g-2}
\hookrightarrow \Sym^{2g-1}$ as an ample divisor and
noticing that the right hand side is a smooth fibre bundle.  It then
follows by the Lefschetz hyperplane theorem.  For the
smooth fibres of $F$, the second statement is an immediate consequence
of the first, whilst for the singular fibres we defer a proof to the
Appendix (\ref{proofforbubbles}).  Finally, the third statement can be
checked by picking a particular complex structure on $\Sigma$ and a
particular rational curve representing $h$.  For instance, if $\Sigma$
is hyperelliptic there is a natural $\pp^1 \subset \Sym^2
(\Sigma)$ which arises from the double covering.  Adding to this a fixed set of
$2g-4$ points of $\Sigma$ defines a rational curve $\phi: \pp^1 \rightarrow
\Sym^{2g-2}(\Sigma)$ whose image is homologous to $h$.  By
construction, the associated closed subset $C_{\phi}$ covers $\Sigma$ and
contains a generic point of $\Sigma$ with multiplicity one.  The
result follows. 
\end{Pf}

\noindent  The Lemma will give us control on
bubbles both for a generic almost complex structure, in the next
Theorem, and for a particular non-generic almost complex structure employed
in the last two sections of the paper.  In any case, with this in
hand, the proof of the following is standard: 

\begin{Thm} \label{Gromovinvariant} Fix some generic $J \in \mathcal{J}$.
There is an integer-valued invariant $\mathcal{I}_r (f)$ which counts
$J$-holomorphic cusp sections of $F:X_r (f) \rightarrow \sss^2$ in the
homology class $[\psi_V] \in H_2 (X_r (f);\zz)$.  It is independent of the
choice of generic $J \in \mathcal{J}$ and of (deformation equivalences of) the
symplectic form $\Omega$. 
\end{Thm}

\begin{Pf}[Sketch] 
The key point will be that for \emph{generic} $J \in \mathcal{J}$ the
moduli spaces of $J$-holomorphic sections of $F$ will be smooth,
compact manifolds of the correct (virtual) dimension.  

\vspace{0.2cm}

\noindent Begin by fixing some
$J \in \mathcal{J}$ and a symplectic form $\Omega$ which tames $J$.
Recall that the defining condition for a map $\phi: \sss^2 \rightarrow
X_r (f)$ to be pseudoholomorphic is that 

\begin{Eqn} \label{J-holequation}
\cdbar_{\phi} \ = \ d\phi \circ j \ - \ J \circ d\phi \ = \ 0.
\end{Eqn}

\noindent Here $\phi$ belongs to a
suitable Banach space of sections (those of H{\"o}lder class $C^{k, a}$
for $k \geq 1$ and $a \in (0,1)$ for instance):

$$ \{ \phi \in C^{k,a}(\mathrm{Maps}(\sss^2; X_r (f))) \ |
\ dF \circ d\phi = \id \}.$$ 

\noindent There is a naturally
induced almost-complex structure on $\phi^* T^{vt} (X_r (f))$, where
the superscript denotes the vertical tangent bundle.  (This bundle is
well-defined on the image of the section, since the
condition that $dF \circ d\phi = \id$ implies that 
$\im (\phi)$ is disjoint from the locus of critical values of $F$.)  The
linearisation of (\ref{J-holequation}) defines a section of Banach
bundles

$$D(\cdbar_{\phi}): \ C^{k,a} (\sss^2; \phi^* T^{vt}X_r(f)) \
\longrightarrow \ C^{k-1,a} (\sss^2; \phi^* T^{vt}X_r(f) \otimes
\Lambda^{0,1} T^*(\sss^2)).$$

\noindent In the usual way, this operator differs from an integrable
$\cdbar$-operator by an order zero operator, and is hence Fredholm
with index 

$$\ind(\cdbar_{\phi}) \ = \ 2 \langle c_1(T^{vt} X_r(f); \Omega),
[\psi_V] \rangle + 2 \dim_{\cc} (F^{-1} (\mathrm{point})).$$

\noindent In our situation this index is zero.  A more general
version of the relevant computation is given in \cite{Sequel}~;  the
key point is the existence of a sequence (well-defined over the image
of the section) 

$$ 0 \rightarrow V \rightarrow T^{vt} X_r(f) \rightarrow W^*
\rightarrow L \rightarrow 0$$

\noindent coming from the differential of the projection $\tau$ of the
Abel-Jacobi map.  The cokernel line bundle $L$ can be shown to be
isomorphic to $K_{\pp^1}$; the result then follows by considering
the associated exact sequences in cohomology, exactly as in the obstruction
computations given in full in section 5.
Recall a \emph{cusp section} is any
cusp curve \cite{McD-S:Jhol} for which the principal component is a
section of $F$.  According to Gromov's compactness theorem, and since
the condition that a curve contain a section component is closed, the moduli
space $\mathcal{M}$ of cusp sections holomorphic with respect to $(j,J)$ is a
compact Hausdorff space (though not necessarily smooth).  Note that
for our class of complex structures $\mathcal{J}$ all bubble
components will lie in fibres of $F$; otherwise consider the
restriction of $F$ to the component and obtain a
contradiction to the homological intersection number with a fibre.

\vspace{0.2cm}

\noindent For ``regular'' $J$ the
moduli space of smooth sections will be smooth and of the expected dimension.
Even though we fix the complex structures to 
be integrable over small discs, there is still a dense notion of
regularity as surjectivity of $D(\cdbar_{\phi})$ will hold whenever
$J$ is generic 
on the bundle restricted to the section away from these discs.  See
(\ref{regularspaces}) and, for proofs, \cite{McD-S:Jhol} and
(\cite{Seidel2}, Section 7).  

\vspace{0.2cm}

\noindent We now claim that, again for a generic - and dense - choice
of the almost complex structure $J$, the moduli space will
also be compact, that is there can be no cusp curves (this is a
"fibred monotonicity" property).  We have already remarked that any
bubbles are vertical.  By (\ref{bubblingfacts}) we know that in
fact any bubble is 
homologous to $Nh$ for some $N \geq 1$.  Any cusp curve has a unique
component with 
non-trivial intersection number $1$ with the fibre of $F$.  We claim
this component is in fact a smooth section.  For it is $J$-holomorphic
with respect to a smooth $J$, and hence by elliptic regularity the map
is everywhere smooth.  Hence it is a section away from the singular
fibres.  Near these fibres, the map $\sss^2 \rightarrow X_r(f)
\stackrel{F}{\longrightarrow} \sss^2$ is holomorphic
with respect to an integrable almost complex structure and hence a
branched covering over its image.  Since
it has local degree one, it must be a local diffeomorphism here also,
and hence the differential of the projection is invertible over the
whole image.  Thus, any non-smooth cusp curve gives rise to a smooth
section of $X_r(f)$.  By the second part of (\ref{bubblingfacts}),
we see that the homology class of this section is $[\psi_V - Nh]$.  But now an
index computation shows that the virtual
dimension for the space of such sections is negative, and so for
generic $J$ such moduli spaces will be empty.  The point is that for
sections of $X_r(f)$ which give rise to cycles in $X'$ in the homology
class $A$, the virtual dimension for the space of $J$-holomorphic
sections is $(A ^2 - K_{X'} \cdot A)/2$.  By the third part of
(\ref{bubblingfacts}),  splitting off a bubble
corresponds to changing $A \mapsto A-[\Fibre]$.  Taking $A = K_{X'}$
gives the required negativity in our case.  The details of this
general index result can be found in \cite{Sequel}, and again involve
studying an exact cohomology sequence coming from the differential of
the Abel-Jacobi map.  An
analogous, but harder, argument will be given in the last section.

\vspace{0.2cm}

\noindent Our compact, zero-dimensional moduli space is naturally
oriented, and our Gromov invariant is the signed count of its number
of points.  Picking a point $P \in \sss^2 \backslash \{ \Crit \}$
defines an evaluation map 

$$ev: \mathcal{M} \rightarrow F^{-1} (P) \ \cong \ \Sym^{r} (\Sigma_g)$$

\noindent and hence a map $ev_*: H_0 (\mathcal{M}) \rightarrow \zz$.
Then we define $\mathcal{I}_r (f) = ev_* ([\mathcal{M}])$. That
this is independent of choices of $J$ and deformations of $\Omega$ is
a familiar cobordism argument, for one-parameter families of moduli
spaces, which asserts that in fact the fundamental homology class
$[\mathcal{M}]$ is independent of such choices.  (Since
$J$-holomorphic curves are smooth whenever $J$ is smooth, we also have
independence of the choice of Sobolev spaces used in constructing the
invariant.)
\end{Pf}

\vspace{0.2cm}

\noindent  In fact, in the rest of the paper, we shall predominantly
work with two almost complex structures from $\mathcal{J}$ that are
\emph{not} regular, and we will provide specific arguments in each
case to obtain smoothness and compactness of the relevant moduli
spaces of holomorphic curves.  We shall make use of the following
standard results: proofs can be found in \cite{Salamon:SW} and
\cite{McD-S:Jhol} respectively. 

\begin{Prop} \label{regularspaces}
\begin{itemize}

\item Let $J \in \mathcal{J}$ be such that the moduli space $\mathcal{M}$ is
a compact smooth manifold of positive dimension $r$ whose tangent
space at every point $\phi$ is given by $\ker (\cdbar_{\phi})$.  Then
the Gromov invariant (virtual class $[\mathcal{M}]^{virt}$) is given
by the Poincar{\'e} dual of
the Euler class of the
\emph{obstruction bundle} over $\mathcal{M}$ with fibre $\coker
(\cdbar_{\phi})$. 

\item Let $U \subset X_r(f)$ be any open set.  By a generic
  perturbation of the almost-complex structure $J$ supported in $U$,
  we can ensure that moduli spaces of all pseudoholomorphic curves
  passing through $U$ are regular.  A similar comment applies to
  structures generic on an open set in $P_r(f)$.
\end{itemize}
\end{Prop}

\noindent In our framework in which the invariant is only defined for
``regular'' almost-complex structures, the first statement above gives
a method for computing $\mathcal{I}$ starting with an irregular but
sufficiently well-behaved $J$.  (The proof is a prototype for the
definition of virtual fundamental classes more generally.)  


\section{Computing the invariant}

In this section, by employing the first result of (\ref{regularspaces}),
we show that the Gromov invariant introduced above is
non-vanishing under a suitable linear constraint on the topology of
the given four-manifold $X$.  

\begin{Prop} \label{goodmodulispace}
There is an almost-complex structure $J_V \in \mathcal{J}$ on $X_r (f)$
for which the moduli space of cusp sections in the class $[\psi_V]$ is a
projective space of dimension $N-1$, where $N = [b_+(X) - 1- b_1
(X)]/2$.  Moreover no sections contain any bubbles.
\end{Prop}

\noindent We start with the following.  Recall that any $\cdbar$-operator on a
vector bundle over a Riemann surface induces a holomorphic structure.
(There is an essentially equivalent discussion of the following in
terms of connexions rather than d-bar operators, if the reader prefers.)
Recall also our notation:  $W$ denotes the ``relative dualising
sheaf'' of the fibration, with fibre canonically identified with the
space of holomorphic one-forms on the fibres, whilst $V$ denotes $W
\otimes \mathcal{O}(-2)$.  

\begin{Lem} \label{indextime}
For a generic d-bar operator $\cdbar$ on $V \rightarrow \sss^2$ the space
of holomorphic sections of the vector bundle has dimension $[b_+(X) -
1 -b_1(X)]/2$.  Every non-zero section defines a smooth section of the
projective bundle $\pp(V)$ in the class $[\psi_V]$.
\end{Lem}

\begin{Pf}
The result is almost an application of standard machinery:  for a
family of $\cdbar$-operators down the fibres of a
smooth fibre bundle, the first Chern class of the index bundle can be
computed using the Atiyah-Singer index theorem.  For a holomorphic
family with singular fibres, this is replaced by the
Grothendieck-Riemann-Roch theorem.  In our case, neither quite apply;
however, in \cite{ivanhodge} an excision argument was used to compute
the first Chern class of the relative dualising sheaf\footnote{In
  \cite{ivanhodge} the determinant of this bundle was denoted by
  $\lambda$.}  $W = f_*
\omega_{X' / \sss^2}$.  From the formulae

$$c_1 (V) = c_1 (f_* \omega_{X' / \sss^2} \otimes \mathcal{O}(-2));
\qquad \rk(V) = g$$

\noindent and, from \cite{ivanhodge}

$$c_1 (f_* \omega_{X' / \sss^2}) \ = \ [\sigma(X') + \sharp 
\{ \mathrm{Crit}(f) \}]/4 $$

\noindent it follows that $c_1(V) = [\sigma(X') + \delta]/4 - 2g$,
writing $\delta$ for the number of singular fibres.  On the other
hand, the Euler characteristic $e(X') = 4-4g+\delta$, and hence
$c_1(V) = [\sigma(X') + e(X')]/4 + (g-1) -2g$.  Since the index is the
sum of the rank and first Chern class, and since the sum of signature
and Euler characteristic for a four-manifold is invariant under
blowing up and down, we arrive at the formula (which is independent of
the degree $k$ of the original Lefschetz pencil)~:

$$\ind (\cdbar) \ = \ \frac{\sigma + e}{4} - 1 \ = \ [b_+(X) - 1 - b_1(X)]/2.$$

\noindent Now on the space of connexions $\{ \cdbar \}$ on a complex vector
 bundle, the discrete-valued function $h^1 (V; J_{\cdbar})$ is
 upper semicontinuous and jumps on a subvariety of positive
 codimension. Hence for generic $\cdbar$ 
 we have $H^1 = 0$ and $H^0 = \ind$; this is a stability result, in
 the lines of Atiyah and Bott's paper \cite{Atiyah-Bott}.

\vspace{0.2cm}

\noindent According to Grothendieck, every holomorphic vector bundle
over $\pp^1$ is a direct sum of line bundles.  Hence $\cdbar$ induces
a splitting of $V$.  By the first assumption in (\ref{highdegree}), we know
that the rank of $V$ is larger than the index~; we deduce $c_1 (V) <
0$.  In this case the 
generic (most stable) splitting is given by 

$$(V, J_{\cdbar}) \ \cong \ \bigoplus_i \mathcal{O}(n_i)$$

\noindent with $n_i \in \{ 0, -1 \}$ for every $i$.  For such a
d-bar operator, the holomorphic sections of $V$ are exactly the constant
sections of the $\mathcal{O}$-factors, and hence \emph{no
  non-zero sections have any zeroes}.  It follows
that every non-trivial section of $V$ defines a section of the
projective bundle $\pp(V)$, in the fixed homology class $[\psi_V]$, which
has no bubbles, since these arise from isolated zeroes of sections of
$V$.  (This is clear, for instance, from a computation in local
co-ordinates near the zero.) 
\end{Pf}

\noindent The discrepancy between the linear constraint $b_+ > 1+b_1$ of
(\ref{maintheorem}) and $b_+ >1$, as required by Taubes \cite{Taubes}, arises
here.  There is a distinguished trivial topological subbundle of $V$
coming from the first homology of the four-manifold $X$ when $b_1 (X)
\neq 0$, and replacing generic d-bar-operators on $V$ by ones
adapted to this subbundle should give the sharper result. We will
suppress the issue in this paper.

\vspace{0.2cm}

\noindent Recall that a choice of a d-bar operator, and hence holomorphic
structure, on $V$ induces one on the manifold $\underline{M}(V)$.
This in turn defines an integrable complex structure on a neighbourhood of
the image of the embedding $\pp(V) \hookrightarrow X_r(f)$, and also
on $W^*$ and hence a neighbourhood of the zero-section (defined by the
canonical lines on each fibre) of the Picard fibration $P_r(f)$.

\begin{Defn}  \label{Jstandardnearzero}
Almost-complex structures $J \in \mathcal{J}$ on $X_r(f)$ and $j$ on
  $P_r(f)$ are \emph{standard
  near the zero-sections} if they agree with the integrable complex
  structures on neighbourhoods of $\pp(V)$ and its image, induced as
  above from a choice of generic d-bar operator on $V$.
\end{Defn}

\noindent In our discussion of almost-complex structures on fibre
  bundles, we noted the freedom to extend structures which
  are prescribed over neighbourhoods of fixed closed sets.

\begin{Lem} \label{allprojectwell}
Given an almost-complex structure $j$ on $P_r(f)$ which is standard 
near the zero section, there is an almost-complex structure $J \in
\mathcal{J}$ on $X_r(f)$ which is  
standard near the zero section and for which the map $X_r(f)
\rightarrow P_r(f)$ is $(J,j)$-holomorphic.
\end{Lem}

\noindent For away from $\pp(V)$ the map is a fibration and we
can just use the existence assertion before (\ref{almostcomplex}).

\begin{Lem} \label{nowheretoproject}
There is an almost-complex structure on $P_r(f)$ which is standard
near the zero-section, and for which any holomorphic section
homologous to the zero-section is itself zero.
\end{Lem}

\begin{Pf}
Along the zero-section, the vertical tangent bundle to $P_r(f)$ is, by
(\ref{itallextends}) 
and the constructions of the Appendix, canonically isomorphic to the
bundle $W^*$.  Recalling that $W = V \otimes K_{\pp^1}$ we have an
identity

\begin{Eqn} \label{indexgames}
\ind_{\cdbar} (V) \ = \ -\ind_{\cdbar}(W^*)
\end{Eqn}

\noindent from Serre duality for any fixed d-bar operators on $V$ and $W$.
If we choose a generic d-bar operator on $V$ as above, and induce the
almost-complex structure on $P_r(f)$ near the zero-section as in
(\ref{Jstandardnearzero}), it follows that $W^*$ has no holomorphic
sections at all, and hence $P_r(f)$ can have no non-zero sections
with image near the zero-section.  But then for a generic extension of
the almost-complex structure from a neighbourhood of the zero-section,
moduli spaces of holomorphic curves in $P_r(f)$ away from zero will be
regular by (\ref{regularspaces}).  Now (\ref{indexgames}) shows that
the virtual dimensions are negative provided $X$ satisfies $b_+ > 1 +
b_1$, and hence the moduli spaces of sections of $P_r(f)$ away from
zero will be empty.
\end{Pf}

\noindent If we combine the three results (\ref{allprojectwell},
\ref{nowheretoproject} and \ref{indextime}) we have a proof of the
proposition (\ref{goodmodulispace}) stated at the beginning of the
section.  That is, choose almost complex structures $j$ on $P_r(f)$
and $J \in \mathcal{J}$ on $X_r(f)$ satisfying
(\ref{nowheretoproject}) and (\ref{allprojectwell}) respectively.  Any
$J$-holomorphic section of $X_r(f)$ projects to a
$j$-holomorphic section of $P_r(f)$, which must be the zero-section.  Hence the
moduli space of smooth $J$-holomorphic sections is just the moduli
space of smooth $J$-holomorphic sections of $\pp(V)$, the
preimage of the zero-section of $P_r(f)$ under the Abel-Jacobi map.
Now by (\ref{indextime}) this last moduli space of holomorphic
sections is a projective 
space $\pp^N$~; this is already compact, and the space
$\overline{\mathcal{M}}$ of $J$-holomorphic cusp sections of $X_r(f)$
is the same projective space.  This is just the required result.

\vspace{0.2cm}

\noindent Appealing again to (\ref{regularspaces}), to compute the
Gromov invariant for sections of $X_r(f)$ in the distinguished class
$[\psi_V]$ we must compute the Euler class of the obstruction bundle
$\mathcal{E} \rightarrow \pp^{N-1}$ (where $N = (b_+ - b_1 -1)/2$).
All the points of the moduli space correspond to sections with image
inside the neighbourhood of $\pp(V) \subset X_r(f)$ which we have
identified with a universal, holomorphic local model
$\underline{M}(V)$, and it follows that we can perform the obstruction
computation here.  By definition, the fibre $\mathcal{E}_{\phi}$ at a
section $\phi$ is given by $H^1 (\nu_{\phi})$ where $\nu_{\phi}$ denotes the
normal bundle to the image.  All the sections $\phi$ are simple, with
no bubbles, and the normal bundle can be identified with the vertical
tangent bundle.  As a piece of notation, we refer to the cokernel of the
natural inclusion $\mathcal{O}_{\mathrm{taut}} \rightarrow
\underline{\cc}^{N+1}$ of the tautological line bundle over projective
space $\pp^N$ into the trivial rank $N+1$ bundle as the \emph{quotient
  bundle}.  Note that this has Euler class $(-1)^N$.

\begin{Prop} \label{quotientbundle}
The obstruction bundle $\mathcal{E} \rightarrow \mathcal{M} \cong
\pp(H^0(V))$ is topologically the dual of the quotient bundle over
projective space.  In particular, $\mathcal{I}_r (f) = \pm 1$ for any
symplectic four-manifold $X$ satisfying $b_+ (X) > 1+b_1 (X)$.
\end{Prop}

\begin{Pf}
Recall the definition of the model $\underline{M}(V)$; we take the
subvariety of $W^* \times_{\pi} \pp(V)$, where $\times_{\pi}$ denotes
the fibre product over $\sss^2$, defined by 

$$\{ (\theta, [x]) \ | \ \theta(x) = 0 \}.$$

\noindent Since the vector bundles $V$ and $W$ can be identified up to
scale, this subvariety is well-defined, and has a projection $\tau$ to $W^*$
whose fibre jumps precisely over $0 \in W^*$.  If we take the
differential of this projection, and recall that the fibre over $0$ is
exactly $\pp(V)$, we obtain a sequence (of holomorphic bundles)

$$0 \rightarrow T\pp(V) \rightarrow T^{vt} \underline{M}(V)
\stackrel{\tau}{\longrightarrow} T(W^*) \rightarrow \cok (d\tau)
\rightarrow 0.$$ 

\noindent We restrict this sequence to a sphere $\phi(\sss^2)$ corresponding to
a point $\phi \in \mathcal{M}$.  It splits into two short exact
sequences

$$0 \rightarrow T\pp(V)|_{\im (\phi)} \rightarrow (T^{vt}
\underline{M}(V))|_{\im (\phi)} \rightarrow \im (d\tau)|_{\im(\tau
  \circ \phi)} \rightarrow 0;$$
$$0 \rightarrow \im (d\tau)|_{\im (\tau \circ \phi)} \rightarrow
T(W^*)|_{\im (\phi)} 
\rightarrow \cok (d\tau)|_{\im (\tau \circ \phi)} \rightarrow 0.$$

\noindent Take the long exact sequences in cohomology: from the first,

$$H^1(T \pp(V)) = 0 \rightarrow \mathcal{E}_{\phi} \rightarrow H^1
(\im d\tau) \rightarrow 0$$

\noindent and from the second

$$H^0 (W^*) = 0 \rightarrow H^0 (\cok (d\tau)) \rightarrow
\mathcal{E}_{\phi} \rightarrow H^0(V)^* \rightarrow H^1(\cok (d\tau))
\rightarrow 0.$$

\noindent  The term $H^0 (V)^* = H^1 (W^*)$ arises by Serre duality.  We claim
that over the image $\phi(\sss^2)$, for $\phi \in \mathcal{M}$, the
bundle $\cok (d\tau)$ defines a copy of the canonical bundle of
$\sss^2$; hence $H^0 = 0$ and $H^1 \cong \cc$ is one-dimensional.
To see this, note that the cokernel of $d\tau$ is generated, for each
$t \in \sss^2$, by a vector in $H^1 (\Sigma_t, \mathcal{O}) = W^*$ on which
$\phi(t)$ is non-zero.  If the vector spaces $V$ and $W$ were dual, we could
choose a metric on $V$ and then the non-zero
vector $\phi(t)$ would trivialise the cokernel bundle; by the same
token, the failure of this is measured by the discrepancy $V = W^*
\otimes K_{\sss^2}$.  Our second sequence above now has the form

$$0 \rightarrow \mathcal{E}_{\phi} \rightarrow H^0 (V)^* \rightarrow
L_{\phi} \rightarrow 0$$

\noindent where $\mathcal{M} = \pp (H^0 (V))$.  The remarks above show
that, after choosing a metric, the section $\phi$ defines a generator
of $L_{\phi}$, and this 
line bundle is therefore dual to the tautological line bundle.  Thus
we have the dual of the standard exact sequence, and the obstruction
bundle is the dual of the quotient bundle as claimed.
\end{Pf}

\noindent It may be helpful to point out a translation of the above if
we work not in the model $\underline{M}(V)$ but on the fibration $X_r
(f) \rightarrow P_r (f)$ itself\footnote{Such a perspective is
  important in \cite{Sequel}, where no analogue of $\underline{M}(V)$
  exists.}.  The cokernel of the canonical
projection from $\Sym^r (\Sigma)$ to $\Pic_r (\Sigma)$ at a divisor
$p_1 + \cdots + p_r$ is canonically equal to $H^1 (\Sigma;
\mathcal{O}_{\Sigma} (\oplus p_i))$.  It follows that the cokernel
bundle, for our situation, is canonically the pullback from the
universal curve $\mathcal{C}_g \rightarrow M_g$ of the bundle $R^1
\pi_* K_{\mathcal{C}}$.  Now Grothendieck's relative duality
\cite{Hartshorne} asserts that for a 
holomorphic fibre bundle $\pi:Z \rightarrow B$ and locally free sheaf $\eta$
on $Z$ there is a natural morphism

\begin{Eqn} \label{relativeduality}
\pi_* (\eta^* \otimes \omega_{Z/B}) \ \rightarrow \ (R^1 \pi_*
\eta)^*
\end{Eqn}

\noindent which is moreover an isomorphism for flat
families of curves.  Taking $\eta = K$, and pulling back the identity
of Chern classes induced by (\ref{relativeduality}), again shows
that $\cok(d\tau)$ gives a copy of the canonical bundle of the sphere
for each $\phi \in \mathcal{M}$.  Then one proceeds as above.  


\section{Strata in the symmetric product}

Our goal is now to construct almost-complex structures on $X_r(f)$ for
which the section provided by the above computation of the Gromov
invariant will yield a symplectic surface.  Our first task is to
review certain natural strata inside the symmetric product of a
Riemann surface.  Fix a single complex curve $(\Sigma, j)$ and a partition 
$\pi: r = \sum a_i n_i$ for integers $a_i, n_i \geq 1$.  Moreover
suppose at least one $a_i > 1$ (to obtain a stratum which is not the
entire space).  Inside $\Sym^r
(\Sigma,j)$ there is a diagonal stratum $\chi_{\pi}$ indexed by
$\pi$ comprising the image of the map

$$\Sym^{n_1} (\Sigma, j) \times \cdots \times \Sym^{n_t}(\Sigma,j) \
  \stackrel{p_{\pi}}{\longrightarrow} \ \Sym^r (\Sigma, j)$$

\noindent which takes tuples $(D_1, \ldots, D_t)$ to $\sum a_i D_i$
(viewed additively as an effective divisor of degree $r$ on
$\Sigma$).  An example: writing $10=3+2+2+2+1$ gives rise to the
mapping

$$\Sigma \times \Sym^3(\Sigma) \times \Sigma \rightarrow
\Sym^{10}(\Sigma): \qquad (p, D, q) \mapsto 3p+2D+q.$$

\noindent For each partition $\pi$ of $r$ of the above form we
obtain a stratum of complex codimension $r-\sum n_i$.  The union of all the
  strata, which is just the image of $\Sigma \times
  \Sym^{r-2}(\Sigma)$ under the map $(p,D) \mapsto 2p+D$, is the
  \emph{diagonal locus} inside $\Sym^r (\Sigma)$.  The various partitions
give rise to a stratification of the diagonal whose top open stratum
comprises the locus where precisely two and no more points of the
tuple have coalesced, and whose open strata parametrise
effective divisors with prescibed multiplicities of points (their
closures allowing multiplicities greater than or equal to those
values).  Although these strata, as closed topological subsets of the
  topological symmetric product, are independent of the complex
  structure $j$, they are not smoothly embedded.  However, each
  stratum $\chi_{\pi}$ has a \emph{smooth model} $Y_{\pi}$,
  which is just the
  domain of the mapping $p_{\pi}$ defined above.  This map is
  finite and holomorphic, and defines an isomorphism over the dense
  open subset of points $(D_1, \ldots, D_t) \in \prod_{i=1}^t
  \Sym^{j_i}(\Sigma)$ for which the supports of the $D_i$ are all disjoint.
  The above discussion is clearly compatible with smooth variations of
  complex structure in the fibres, so at least over $\sss^2 \backslash
  \{ \Crit \}$ we have smooth fibre spaces $Y_{\pi} \rightarrow
  \sss^2_*$ together with maps, which are fibrewise holomorphic, to the strata
  $\chi_{\pi} \rightarrow \sss^2_*$.  (The situation is rather
  similar to the case of the map $t \mapsto (t^2, t^3)$ which takes
  $\cc$ homeomorphically onto a singular complex curve in $\cc^2$.)

\vspace{0.2cm}

\noindent Along with the diagonal strata, we need to consider strata
arising from the exceptional spheres of the fibration $f: X'
\rightarrow \sss^2$.  Here, for a fixed point $e \in \Sigma$ and
multiplicity $a(e) \in \zz_{>0}$, we have
a stratum $\chi_{a(e)} \subset \Sym^r (\Sigma)$ defined as the image of the map

$$\Sym^{r-a(e)}(\Sigma) \ \stackrel{p_{a(e)}}{\longrightarrow} \ \Sym^r
(\Sigma): \qquad D \mapsto D + a(e)e.$$

\noindent The image of the map is again not a smooth subset of the
$r$-th symmetric product, but we again regard the domain of $p_{a(e)}$ as a
smooth model for $\chi_{a(e)}$, and again $p_{a(e)}$ is holomorphic, a
homeomorphism on a dense open set, and compatible with passing to
smooth families.  Combining the two discussions, we define strata
$\chi_{\pi, \aleph} \subset \Sym^r(\Sigma)$ as we vary over
partitions $\pi$ of the 
integer $r-|\aleph|$ and ordered subsets $\aleph
\subset \{ 1, \ldots, k^2 \omega_X^2 \}$ of the set of exceptional
spheres (which may be taken with multiplicities, so $\aleph$ may have
repeated elements).  There are smooth models

$$Y_{\pi, \aleph} \rightarrow \chi_{\pi, \aleph}$$

\noindent where if $\pi = \sum_{i=1}^t a_i n_i = r-|\aleph|$ we map 

$$\prod_i \Sym^{n_i} (\Sigma) \rightarrow \Sym^r (\Sigma); \qquad
(D_1, \ldots, D_t) \mapsto \sum a_i D_i + \sum_{j \in \aleph} e_j.$$

\noindent Here $e_j$ denotes the point of intersection of the
exceptional sphere $E_j \subset X'$ with the fibre corresponding to
$\Sigma$.  The complex codimension of the stratum is $r-\sum n_i -
|\aleph|$ where the elements of $\aleph$ must be counted with their
multiplicities.  As we vary in families, we obtain smooth fibre bundles
$Y_{\pi, \aleph}$ over $\sss^2_*$.  The importance of these for
the proof of (\ref{maintheorem}) is explained by the following remarks
(which rely on a lemma at the start of the next section but may better
motivate the present discussion than be deferred).

\vspace{0.2cm}

\noindent  Suppose we have a  smooth section  $\phi$ 
of $X_r(f)$.  There is then an obvious way to associate a 
closed set $C_{\phi} \subset X'$ to $\phi$.  Note first that this
cannot meet the  
singular set in the special fibres. Suppose for simplicity that $\phi$
also lies in the dense 
open set $\Sym^r(\Sigma_0 \backslash \{Q \})$ at each of the
singular fibres, where $Q$ denotes the node, and that $\phi$ is
disjoint from all the strata $\chi$ of real codimension greater than
two.  Clearly $C_{\phi}$ is a 
smooth embedded surface away from the points where $\phi$ meets the 
diagonals. 

\vspace{0.2cm}

\noindent Furthermore, we claim that $C_{\phi}$ is a standard 
surface (\ref{standardsurface}) if $\phi$ meets the (real) codimension
$2$-strata transversally,  
with locally positive intersection.  (Away from the strata
of higher codimension, the diagonals are smoothly embedded, and so
transversality here makes good sense.)  The only real codimension two
strata correspond to the partition $r = 2+1+\cdots+1$ and to the
multiplicity one loci $\chi_{a(E)=1}$.  In each case, by throwing away
local smooth sections, we can reduce
to a model inside the second symmetric product.  For the first case,
after a diffeomorphism the model is as follows.  The curve is given by 

$$\{ x^2 + y = 0 \} \ \subset \ \cc^2 \ \rightarrow \ \cc;$$

\noindent the map to $\cc$ is projection on the second factor, and the
associated map from the base to the second symmetric product of some
fixed fibre $F$ is  given by $t \mapsto [i\sqrt{t}, -i\sqrt{t}]_{\sim}$.
(The equivalence relation $\sim$ is that of the natural covering $F \times F
\rightarrow \Sym^2(F)$.)  The image of the map to $\Sym^2(F)$ has
intersection number $1$ with the diagonal, and hence since the model
is holomorphic this must be a transverse intersection.  By contrast,
if for instance $C_{\phi}$
contains a nodal point we can
take a neighbourhood of the node~:

$$\{ x^2 + y^2 =0 \} \ \subset \ \cc^2 \ \rightarrow \ \cc$$

\noindent with the fibration again being second projection.  Then the map
from a disc to the second symmetric product of the fibre $F$ is locally
\emph{two-to-one};  $\pm t \mapsto [it, -it]_{\sim}$.  Now the
intersection with the diagonal has multiplicity two, and amounts to a
smooth tangency.   We can always perturb this second
model to two tranverse intersections with the diagonal;
for instance, $\{ x=0 \} \subset \cc^2$
passes through the node of a local map $(x,y) \mapsto x^2 + y^2$ but
the perturbation $\{ x = \varepsilon \}$ is instead tangent to the
fibres at two points.  Similarly, intersections of $\phi$ with the
stratum $\chi_{a(E)=1}$  at a point $F(e) \in \sss^2$ correspond to
transverse intersections  of the
surface $C_{\phi}$ with the fixed section $E \ni e$.  Thus if the
section is transverse to all the 
strata of the diagonal it defines a smooth symplectic curve which
meets $E$ locally positively.  More generally we have

\begin{Prop} \label{buildstandardsurface}
Suppose $\phi$ lies inside a stratum $\chi_{\pi, \aleph}$, does not meet any 
stratum of (real) codimension bigger than $2$ in $\chi_{\pi,
  \aleph}$ and cuts the codimension $2$ strata transversally and 
with locally posistive intersections. Then $C_{\phi}$ is a union of standard 
surfaces in $X$, with positive transverse intersections and no triple
intersections.  Moreover, at an intersection point, either one of the
components must be an exceptional sphere, or the two components must
have strictly different multiplicities.
\end{Prop}

\noindent It is not quite clear how to interpret transversality here
(consider, again, the case of a disc in $\cc^2$ lying inside the
cuspidal cubic curve $\{ y^2 = x^3 \}$ and the meaning of
transversality of the disc to the origin).  However, as we shall explain
at the start of the next section, $\phi$ defines an associated
section $\tilde{\phi}$ of a unique minimal smooth model $Y_{\pi,
  \aleph}$ in which the hypotheses assert that it is disjoint from the
strata of real codimension at least four.  In this case,
transversality makes sense for the smooth open loci in the codimension
two strata.  (This amounts to pulling a disc back to the smooth domain
$\cc$ of the cuspidal cubic in our analogy, and asking for
transversality to the origin here.)  Clearly
$C_{\tilde{\phi}}$ is naturally partitioned by collecting together the
components 
belonging to each disctinct factor
in the partition $\pi$ which indexes its associated stratum.  This 
exactly groups the components by their multiplicities.  The
discussion before the Proposition implies that
the union of the components of any given multiplicity is a smooth
symplectic surface, and moreover that the
components coming from different factors can
meet only transversally.  By construction
$C_{\phi}$ is given, as a set, by the components of
$C_{\tilde{\phi}}$ and various copies of exceptional curves.  The
result follows.

\begin{Rmk}
It is worth noting that although the two kinds of strata $\chi_{\pi}$ and
$\chi_{a(E)}$ play similar roles in the proof, conceptually they enter
with different flavours~: we are trying to find sections
which will be transverse to all the strata of the first sort, but
which for
topological reasons will necessarily lie in all the strata
$\chi_{a(E)=1}$.  Many of the subsequent technical difficulties arise
  because we cannot exclude the possibility that in fact the sections
  \emph{do} also lie inside various of the $\chi_{\pi}$ or
  $\chi_{a(E)>1}$.
\end{Rmk}

\vspace{0.2cm}

\noindent  We can also attach a multiplicity $a_{i}$ to each component
$C_{i}$ of $C_{\phi}$, derived from its multiplicity in the symmetric
product (so for the exceptional sections, these multiplicities are the
$a(E)$ we have encountered before).
This associates to $\phi$ a positive symplectic divisor in the sense
of (\ref{possympldivisor}). 
We associate the homology class $[C_{\phi}]=\sum a_{i}[C_{i}]\in 
H_{2}(X')$ to $\phi$, and if $\phi$ is homotopic to a section of
$\pp(V) \subset X_r(f)$ in the homotopy class $[\psi_V]$ of sections
coming from $V$, then $[C_{\phi}] = K_{X'}$. 

\noindent The smooth models $Y_{\pi, \aleph}$ form smooth fibre
bundles as we vary the Riemann surface over $\sss^2_*$, and these
admit almost-complex 
structures for which the projection map to the sphere is
pseudoholomorphic and which 
agree with the standard integrable structures on the fibres.  The
proof of this precedes as in (\ref{almostcomplex}), again building on
results of the Appendix.  

\begin{Defn} \label{compatiblewithstrata}
An almost-complex structure $J \in \mathcal{J}$ on $X_r(f) \rightarrow
\sss^2_*$ is 
\emph{compatible with the strata} if there are almost-complex
structures on all the $Y_{\pi, \aleph}$ as above for which the
canonical projection maps $Y_{\pi, \aleph} \rightarrow X_r(f)$
are holomorphic.
\end{Defn}

\noindent  Again, such a $J$ will tame the symplectic forms $\Omega_R$
for $R> R(J)$.  Then the important ingredient for us is the 

\begin{Prop} \label{Jcompatiblewithstrata}
There exist almost-complex structures on $X^*_r(f)$ which are
compatible with the strata and which agree with the given (integrable)
structures near the ends.
\end{Prop}

\noindent We shall give the proof for a stratum $\chi_{\pi}$ with
$|\aleph| = 0$; the general case is complicated only by notation.

\vspace{0.2cm}

\begin{Pf}
Recall that for an open set $U$ in $\cc$ we have the map
$\usigma: \Sym^{r}(U)\rightarrow \cc^{r}$ defined by the elementary symmetric 
functions.   Let $\pi = \sum_{i=1}^t a_i n_i$ be a partition  of
$r$; write $\sum_{i=1}^t n_i = s$. We then have sequences of maps 

$$Y_{\pi}(U) \rightarrow \usigma_{\pi}(U) \subset \cc^s$$
$$Y_{\pi}(U) \rightarrow \Sym^{r}(U)\rightarrow \usigma(U)\subset
\cc^{r}. $$ 

\noindent Here the first map is defined by using the identification
with the appropriate product of smaller 
symmetric products.  The key point is that if we have a holomorphic 
diffeomorphism $\alpha:U \rightarrow V \subset \cc$ we get a diagram of 
holomorphic maps connecting the above $U$-sequences and $V$-sequences
(where all the vertical maps are induced by $\alpha$):

$$
\begin{array}{ccccccccccc}

\cc^s & \supset & \usigma_{\pi}(U) & \leftarrow & Y_{\pi}(U) &
\rightarrow & \Sym^r
(U) & \rightarrow & \usigma(U) & \subset & \cc^r \\
& & \downarrow & & \downarrow & & \downarrow & & \downarrow & & \\
\cc^s & \supset & \usigma_{\pi}(V) & \leftarrow & Y_{\pi}(V) &
\rightarrow & \Sym^r (V) & \rightarrow & \usigma(V) &
\subset & \cc^r 

\end{array}
$$

\noindent Here the vertical map on the left is just the product of the
maps induced on  
the smaller symmetric products. Now apply this to the charts on
$X^*_r(f)$ defined by the restricted charts on $X^{*}$.
These charts map sets   
$ G\times D$ into $X^*_r(f)$, where $G$ is a product of open sets 
of the form $\usigma(U_{\alpha})$, and $U_{\alpha}$ are local
holomorphic charts  
along the fibres of $X^{*}$.   In one of these charts, say $G_{1}\times D$,
we have a preferred 
almost-complex structure $J_{1}$ given by the product structure.
Let $G_{2}\times D$ be
another such chart, with the same image. These charts compare by a smooth
family of holomorphic maps $g_{\tau}:G_{1}\rightarrow G_{2}$. Now these maps
$g_{\tau}$ are not arbitrary holomorphic maps; they are all induced in the 
manner above by  holomorphic diffeomorphisms between open sets in $\cc$. 
Hence these maps fit into a diagram:

$$
\begin{array}{ccc}
G_{1,\pi}\times D & \rightarrow & G_{1}\times D \\
\downarrow & & \downarrow \\
G_{2,\pi}\times D & \rightarrow & G_{2}\times D 
\end{array}
$$

\noindent where $G_{1,\pi}, G_{2,\pi}$ are open sets in some $\cc^s$
which give  
local charts on the smooth model $Y_{\pi}$ of the stratum $\chi_{\pi}$. 
The almost-complex structure $J_{1}$ is represented in the second
chart $G_2 \times D$ by a 
vector field along the fibres $\xi= \frac{\partial g_{\tau}}{\partial 
\overline{\tau}}$. The product structure in the first chart also gives
a local 
almost-complex structure $J_{1,\pi}$ on the smooth model of the stratum. Thus 
there is another vector field $\xi_{\pi}$ along the fibres of $G_{2,\pi}$ and 
clearly the derivative of the map $\mu: G_{2, \pi} \times D \rightarrow G_2
\times D$ takes $\xi_{\pi}$ to $\xi$.

\vspace{0.2cm}

\noindent  Now for any smooth function $\beta$ on $G_{2}\times D$, the 
product $\beta \xi$ defines a vector field and hence another
almost-complex structure on 
$G_{2}\times D$. The composite of $\beta$ with the map $\mu$  
is a smooth function
$\beta_{\pi}$ on $G_{2,\pi}\times D$.
So $\beta\xi$ and $\beta_{\pi}\xi_{\pi}$ define local almost-complex 
structures on $X^*_r(f)$ and $Y_{\pi}$, respectively,
which are compatible, in the sense that the natural projection map is
holomorphic. 

\vspace{0.2cm}

\noindent With these remarks in place we can carry through the usual
procedure to construct a compatible almost-complex structure on the whole of 
$X^*_r(f)$: we take the trivial structures in charts
given by the 
product structure and glue these together using a partition of
unity. Then the key observation is merely that
a smooth partition of unity on
$X^*_r(f)$ induces, by pullback, a smooth partition of
unity on the smooth models of the strata.
We can arrange that the complex structure agrees with the standard ones 
near the ends of $\sss^2_*$ by
taking the original restricted charts on $X'$ to be holomorphic near
the singular fibres.
\end{Pf}

\noindent We now need to discuss how to extend these strata over the
singular fibres.  From one perspective, the diagonal strata arise from
the semigroup structure on $\amalg_n \Sym^n (\Sigma)$, for which there is
no analogue for the union of the zero-fibres of the $X_n(f)$.  This
union, as observed in the Appendix, is just $\amalg_n \Hilb^{[n]}
(\Sigma_0)$, where $\Sigma_0$ is the fibre of $f$ over the critical
value.  On the
other hand, we do have a canonical projection map $\Hilb^{[n]}(Z)
\rightarrow \Sym^n (Z)$ for any space $Z$ \cite{Nakajima}, and using
this we can pull back the strata defined by the semigroup action on the
singular symmetric product spaces.  In the usual way, the diagonal
strata will be (singular) algebraic varieties in the total space of
$\mathcal{X}_r(f)$ given a holomorphic projection from a smooth
complex surface $\mathcal{X}$ to the disc $D$ as in
(\ref{itallextends}).  The point of importance for us is the following

\begin{Lem} \label{stillcodim2}
For every partition $\pi$ of $r$, the stratum $\chi_{\pi}$ of the zero-fibre
$(\mathcal{X}_r(f))_0$ meets $\Sym^r (\Sigma_0 \backslash \{Q\})$ in a
Zariski dense set.  In particular, the complement of this set has
codimension at least two in every stratum.
\end{Lem}

\noindent This follows, for instance, from the discussion in
Nakajima's notes (\cite{Nakajima} or the Appendix) of the
Hilbert-Chow morphism, and the fact that the analogous statement for
the strata in $\Sym^r$ does hold.


\section{Constructing the symplectic surface}

In this section we assemble the pieces already established to complete
the proof.  The delicacy will be that although we consider sections
which have
index zero in $X_r(f)$, their indices as sections of the strata $Y_{\pi,
  \aleph}$ (in which \emph{a priori} they may lie) are not known.  To
get around this, we will establish a regularity result for almost
complex structures compatible with the strata which will have the
following consequence.  Fix some $J$ compatible with the strata and a
holomorphic section $\phi$ provided by the non-vanishing of the Gromov
invariant.  (Suppose there are no bubbles, a fact we will prove later
in the section.)  If the section lies inside some stratum $\chi$ we pass
to an associated smooth section of a smooth model $Y$.  If this section
has negative index, our regularity result will allow us to perturb $J$
and assume no such section in fact existed.  Hence, we move $\phi$
outside of the stratum $\chi$. If, on the other hand, the index for
the associated section inside $Y$ is non-negative, our regularity result
on $Y$ will allow us to assume that $\phi$ is
transverse to all the smaller strata.  An obvious finite induction
(successively pushing off any strata of negative index), coupled with
the fact that an intersection of dense 
sets is dense, will enable us to conclude there are "enough" almost complex
structures $J$ compatible with the strata~; then, as usual, Sard's theorem
provides such for which smooth holomorphic
sections satisfy all the conditions of proposition
(\ref{buildstandardsurface}).  

\vspace{0.2cm}

\noindent Before turning to this programme, let us review more
carefully the way we associate the sections of the smooth models.
We begin with a discussion of
pseudoholomorphic maps to stratified spaces.  
Let $B$ be an open ball in $\cc^{n}$ and $A$ be a complex analytic 
subvariety in $B$. Thus $A$ has a stratification by subsets 
which are locally-closed complex manifolds in $B$.
Suppose $\mu$ is a smooth almost-complex structure on
$B$. We say that $A$ is a $\mu$-subvariety if all the strata of $A$ are
complex submanifolds with respect to the almost-complex structure
$\mu$.  Equip the unit complex disc $D$ with its standard complex structure.

\begin{Lem}
Let $A\subset (B, \mu)$ be a $\mu$-subvariety in the above sense. 
If $f:D\rightarrow B$ is pseudoholomorphic with respect to the structure 
$\mu$ then the set $f^{-1}(A)$ is either the whole of $D$ or a discrete subset 
of $D$.
\end{Lem}

\begin{Pf}
By considering the stratification of a singular variety we can reduce to the 
case when $A$ is a submanifold. By applying a suitable holomorphic 
diffeomorphism we can then suppose 
that $A$ is a linear 
subspace in $\cc^{n}$. In the case when $A$ is a complex line the result 
is proven by McDuff in (\cite{AudLaf:eds}, Chapter 6). For the general
case one finds that the
proof of \cite{AudLaf:eds} goes over unchanged. 
\end{Pf}

\noindent An obvious extension of the above gives

\begin{Cor}
Suppose $A_{1},A_{2}$ are $\mu$-subvarieties in $B$ as above, neither of which 
is contained in the other. Then if $f:D\rightarrow B$ is a pseudoholomorphic 
map whose image lies in $A_{1}\cup A_{2}$ then either $f$ maps into $A_{1}\cap 
A_{2}$ or $f$ maps into precisely one of the sets $A_{1}$ and $A_{2}$.
\end{Cor}

\noindent We can apply these results to our subsets
$\chi_{\pi, \aleph}$,  
since these are represented in local charts by complex varieties. We
deduce that 
there is a unique minimal stratum $\chi_{\pi, \aleph}$ associated to
a pseudoholomorphic 
section $\phi$ of $X_r(f)$. That is, $\phi$ lies in
$\chi_{\pi, \aleph}$ but not in any smaller stratum.  Moreover
$\phi$ meets the  
smaller strata in discrete (hence finite) sets.
 
\begin{Lem}
If $\chi$ is the ``minimal'' stratum associated to a section $\phi$ as 
above, then there is a unique holomorphic section $\tilde{\phi}$ of
$Y_{\chi}$ which maps to $\phi$ under the canonical map $Y_{\chi}
\rightarrow X_r(f)$.
\end{Lem}

\noindent The map $Y_{\chi} \rightarrow X_r(f)$ is a homeomorphism on
a dense open set, so $\tilde{\phi}$ is  
uniquely defined as a continuous section on a dense set.  Moreover, it
is smooth  and pseudoholomorphic away from the lower strata, i.e. where the
derivative of the projection of $Y_{\pi}$ is injective.  In a local
chart near one of the finitely many intersection points with the lower
strata, $\tilde{\phi}$ is bounded and hence
extends to a smooth map on the entire sphere.  Thus 
the assertion follows from the fact that a continuous map which satisfies a 
pseudoholomorphic mapping equation  outside a discrete set is actually 
pseudoholomorphic everywhere.  This follows from elliptic regularity
as in \cite{Sacks-Uhl}, for example (since we know $\tilde{\phi}$ is
\emph{a priori} continuous we are in the easy case).

\vspace{0.2cm}

\noindent In order to achieve the transversality of
(\ref{buildstandardsurface}) we will establish a regularity result for
almost-complex structures which are compatible with the strata.  Let
$\mathcal{S} \subset \mathcal{J}$ denote the set of such
almost-complex structures on 
$X_r(f)$;  this is an affine space in the familiar way.  For each $j
\in \mathcal{S}$ we have an almost-complex structure $j_{\pi, \aleph}$ on
$Y_{\pi, \aleph}$ for which the projection $Y_{\pi, \aleph}
\rightarrow X_r(f)$ is holomorphic. 
Let $H_{\pi,\aleph}$ be the set of homotopy classes of sections of
$Y_{\pi,\aleph}$, so for each 
$h\in H$ and $j\in \mathcal{S}$ we have a moduli space $M_{h,j}$ of
pseudoholomorphic sections of $Y_{\pi, \aleph}$ which do not all lie in
the closure of any proper stratum in $Y_{\pi, \aleph}$. 
Let $M'_{h,j}$ denote the subset of sections 
which are ``good'', that is which are 
(i) transverse to all the lower strata in $Y_{\pi, \aleph}$ and
(ii) lie in the dense open set $p_{\pi, \aleph}( \prod_i \Sym^{n_i}
(\Sigma_0 \backslash \{ Q \}))$ over each critical value of $f$.  The
vanishing of the index for our original problem does not allow us to
conclude anything about the dimensions of these spaces $M_{h,j}$.  However~:

\begin{Prop} \label{Jregularforstrata}
For generic $j \in \mathcal{S}$ all pseudoholomorphic sections of
$(Y_{\pi, \aleph},j_{\pi, \aleph})$  
are regular, so $M_{h,j}$ is a manifold of the expected dimension.
Moreover, for generic $j$ and all $h$, the set $M'_{h,j}$ is  dense in
$M_{h,j}$.
\end{Prop}

\noindent The corresponding assertion where $j$ varies in the set of
\emph{all} almost-complex structures on $Y_{\pi, \aleph}$ (or even all
those structures compatible with its fibration) is
standard. The point here is that we are 
only allowed to consider the restricted set of almost-complex
structures which arise from the 
compatible structures on $X_r(f)$.  As usual, the proof proceeds by
constructing a smooth "universal" moduli space~; then the appropriate
class of generic $j$ are given by the regular values of a projection
map to the space 
$\mathcal{S}$.   We will ignore the Sobolev spaces, which are
standard, and focus on the key geometric feature of the argument.  Fix
some reference structure  
$j^{(0)} \in \mathcal{S}$, inducing an almost-complex structure
$j^{(0)}_{\pi, \aleph}$ on $Y_{\pi, \aleph}$.  Recall that any other
almost-complex  structure  on $Y_{\pi, \aleph}$ differs by a vertical vector
field  on  $Y_{\pi, \aleph}$, as in the discussion of section 4.  Let
$\phi$ be a $j^{(0)}_{\pi, \aleph}$-holomorphic section of 
$Y=Y_{\pi,\aleph}$.

\vspace{0.2cm}

\noindent  Let $\cdbar$ denote the operator $(s,J) \mapsto
\cdbar_J (s)$ viewed as a map of Banach manifolds.
For the universal moduli space of holomorphic curves to be
smooth, we need to know 
that the cokernel of the linearised operator $D(\cdbar)$ 
is everywhere zero.  If there is some non-zero
element of such a cokernel, then there is a non-zero element $\eta \in
(\Lambda ^{0,1} T^* \sss^2) \otimes \phi ^* TY$ of the kernel of
the adjoint map.  By the Hahn-Banach theorem, this is possible only if
$\eta$ is orthogonal to the image of $D(\cdbar)$.  Clearly it
is enough to prove that in fact all such $\eta$ vanish on a dense set, and so
we can assume for contradiction that $\eta$ is non-vanishing at $q \in
\Delta_{\phi}$, where $\Delta_{\phi}$ is the
open set in $\sss^2$ obtained 
by removing the critical values of $f$ and the points $p$ where
$\phi(p)$ lies in some smaller stratum of $Y$.  
Following the argument of
(\cite{McD-S:Jhol}, p.35) the crucial point is to construct a tangent
vector $v \in \mathrm{End}(TY,j^{(0)}_{\pi, \aleph})$ to the space of
almost complex structures for which 
$\langle \eta, v \circ d\phi \circ j_{\sss^2} \rangle$ is
non-vanishing.  (The RHS of the inner product is one term in
$D(\cdbar)$ evaluated at a tangent vector of the shape $(\delta \phi,
v)$, and the non-vanishing contradicts the orthogonality to the image
of $D(\cdbar)$.)  Hence it is certainly enough to show that \emph{all}
tangent vectors 
to $(T_q Y, (j^{(0)}_{\pi, \aleph})_q)$, viewed as a point of the space of
almost complex structures on $Y$ at $q$ arising from structures compatible
with the fibration, can be
generated by perturbations inside the space $\mathcal{S}$.  This is
the content of the following~:

\begin{Lem}
There is an open neighbourhood $\Delta \subset \sss^{2}$ of $q$ such that 
for any compactly supported section $v$ of $\phi^{*}(T^{vt}Y)$ 
over $\Delta$, there is a $j\in \mathcal{S}$ such that
$j_{\pi,\aleph}$ differs from $j^{(0)}_{\pi,\aleph}$ by a vertical
vector field $\xi$ on $Y$ which pulls back to $v$ over $\Delta$. 
\end{Lem}

\begin{Pf}
We return to the charts for $X_r(f)$ and $Y$ obtained 
from restricted charts on $X'$. We have seen that a holomorphic 
diffeomorphism $\alpha:U\rightarrow V$ of open sets in $\cc$ induces 
holomorphic diffeomorphisms from $\usigma(U)$ to $\usigma(V)$ and
from $\usigma_{\pi}(U))$ to $\usigma_{\pi}(V)$. 
In the same way a holomorphic vector field $\eta$ on $U$ induces  holomorphic 
vector fields 
$\ueta$ on $\usigma(U)$ and $\ueta_{\pi}$ on
$\usigma_{\pi}(U)$. A smooth family $(\eta^{\tau})$ of holomorphic vector 
fields on $U$, smooth in the base parameter $\tau$, then gives vertical vector 
fields over $\usigma(U)\times D$ and
$\usigma_{\pi}(U)\times D$ which  
we can use to deform the given almost-complex structure
$j^{(0)}$. Suppose we have a point $x$ of
$\usigma_{\pi}(U)$ which does not lie in any lower stratum. Then it is easy to 
see that 
for any 
tangent vector $w$ to $\usigma_{\pi}(U)$ at $x$
there is a holomorphic vector field $ \eta$  on $U$ such that  
$\ueta_{\pi}$ is equal to $w$ at $x$.
When one unwinds the definitions this just amounts to finding a vector field 
on $U$ taking prescribed values at a finite set of points.
More generally, suppose $g$ is a smooth 
map from the disc into $\usigma_{\pi}(U)$ with $g(0)=x$. Then one can find a 
neighbourhood $\Delta$ of $0 \in D$ and  smoothly varying families of vector 
fields $(\eta_1^{\tau}, \ldots, \eta_s^{\tau})$ on $U$, depending on the
parameter 
$\tau\in \Delta$, such that for each $\tau$ the corresponding vector fields
$(\ueta_{1,\pi}^{\tau},\dots,\ueta_{s,\pi}^{\tau})$ on
$\usigma_{\pi}(U)$ evaluated 
at $g(\tau)$ give  a basis for the tangent space of $\usigma_{\pi}(U)$
at $g(\tau)$.  This construction leads immediately to the proof of the
Lemma, when  
we take a co-ordinate chart around the point $\phi(q)$ in $X_r(f)$, and extend
using cut-off functions.
\end{Pf}

\noindent Both statements of (\ref{Jregularforstrata}) follow from the
above.  Using the affine structure on $\mathcal{S}$ induced by vector
fields, we can superimpose these local deformations.  As a
consequence, we have the statement~:

\begin{Lem} \label{pushsections}
Let $\phi$ be a $j^{(0)}_{\pi, \aleph}$-holomorphic section of
$Y=Y_{\pi,\aleph}$.  Then every section of
$Y$ which co-incides with $\phi$ outside a compact subset of $\Delta_{\phi}$
is $j_{\pi, \aleph}$-holomorphic for some $j \in \mathcal{S}$.
\end{Lem}

\noindent The relevance of this result is the following. Fix some
symplectic form $\Omega$ on $X_r(f)$ and choose a compatible $J \in
\mathcal{S}$.  If we are given a $J$-holomorphic 
section $\phi$ then any sufficiently $C^k$-small perturbation $\phi'$ of $\phi$
(which co-incides with $\phi$ outside of $\Delta_{\phi}$) will be
holomorphic for some
perturbed $J' \in \mathcal{S}$ which still tames the fixed $\Omega$.
Using this, it follows that the evaluation map from the space of
$J$-holomorphic sections to a fibre of $X_r(f)$ is submersive~; this
is standard for regular $J$, and important for defining the invariants
coming from higher 
dimensional moduli spaces (some of these are computed in \cite{Sequel}).

\vspace{0.2cm}

\noindent To complete the proof, we must provide a discussion of
bubbling for the almost complex structures compatible with the strata.
This shall hinge entirely on the assumption (\ref{highdegree}) that we
originally fixed a Lefschetz pencil of (sufficiently) high degree $k$.  
Fix a generic almost-complex
structure $J \in \mathcal{S}$ on $X_r(f)$ which is compatible with the
strata.  From the 
non-triviality of the Gromov invariant, we know the moduli space
$\overline{\mathcal{M}}[\psi_V]$ of holomorphic cusp
sections in the homology class $[\psi_V]$ is non-empty.  Fix such a
cusp section $\overline{\phi}$, and denote by $\phi$ its unique
section component (which 
is smooth by the arguments of the previous section).  This section $\phi$
has an associated minimal stratum $Y$, where it defines a point
$\tilde{\phi}$ in a
moduli space $\mathcal{M}_Y$ of pseudoholomorphic sections.  Since $J$ is
generic, there is a dense set in $\mathcal{M}_Y$ corresponding to
sections which are ``good'' over the singular fibres and which meet
all lower strata transversely.  Suppose $(\phi, \tilde{\phi})$ are
good in this sense.

\begin{Lem}
The section $\overline{\phi}$ cannot contain any bubbles, so
$\overline{\phi} = \phi$.
\end{Lem}

\begin{Pf}
We have already established (\ref{bubblingfacts}) that all bubble
components are homologous to multiples of the standard projective line
$h$ inside the projective space $\pp(V)$. It follows that if there are
any bubbles, then the section component $\phi$ of the curve lies inside a
moduli space of curves $\mathcal{M}[\psi_V-Nh]$ in the homology class
$[\psi_V] - Nh$ for some positive integer $N$. 

\vspace{0.2cm}

\noindent The section $\phi$ of $X_r(f)$ defines a section
$\tilde{\phi}$ of some associated minimal stratum $Y$.  We claim that
in fact $C_{\phi}$ must contain all of the exceptional sections of the
Lefschetz fibration, so in particular $Y$ is a smooth model of a
stratum lying inside $\cap_i \chi_{a(E_i) = 1}$.  By symmetry, it is
enough to prove this for some fixed exceptional curve $E$.  If
$C_{\phi}$ does not contain $E$, then $\tilde{\phi}$ is transverse to
$Y \cap \chi_{a(E)=1}$, where the last term denotes the obvious real
codimension two stratum of $Y$.  According to the discussion before
(\ref{buildstandardsurface}), in this case $C_{\phi}$ has locally
positive intersections with the exceptional curve $E$, and in
particular $C_{\phi} \cdot E \geq 0$.  But this is a contradiction~:
we know from (\ref{bubblingfacts}) that $C_{\phi}$ represents the
homology class $PD[K_{X'}] - N[\Fibre]$ in $H_2 (X')$, and this has
intersection number $-(N+1)$ with $E$.

By combining (\ref{buildstandardsurface}) and
(\ref{bubblingfacts}), we see that the section component
$\tilde{\phi}$ of the 
associated minimal stratum defines a smooth symplectic surface, and
after we pass back to $C_{\phi}$ and assign multiplicities we obtain a
positive symplectic divisor on
$X'$ in the homology class $PD[K_{X'}]-N[\Fibre]$.  Moreover this contains all
of the exceptional components to multiplicity at least one.  On the
other hand, if it contains some exceptional curve $E$ to multiplicity
greater than one, then - running the same intersection argument as
before and using (\ref{buildstandardsurface}) again - we see that the
other components meet $E$ transversely and positively.  We can now
apply our long-ignored smoothing lemma (\ref{smoothingtwo}) to
separate out the exceptional curves and obtain a smooth symplectic
surface inside $X'$ which (i) contains the exceptional curves to
multiplicity precisely one and (ii) represents the homology class
$K_{X'} - N[\Fibre]$.  Throwing out the exceptional curves, we can
push the resulting surface down to $X$. This remains symplectic, by 
(\ref{symplecticforms}) and the ensuing discussion~; but now we have a
contradiction, for the surface in $X$ represents the homology class
$K_X - N [(k/2\pi) \omega_X]$.  But by our initial choice of $k$ in
(\ref{highdegree}) the evaluation of $\omega_X$ against this class is
negative, and so the existence of such a symplectic surface is precluded.
\end{Pf}

\noindent This essentially completes the proof.  The non-vanishing of
the Gromov invariant for sections in the class $[\psi_V]$, together
with the regularity results above and the absence of bubbles, means
that we have a smooth holomorphic section $\phi$ for which the
associated positive symplectic divisor $C_{\phi}$ is in the class $K_{X'}$.  By
the same argument as above, this divisor contains all of the
exceptional sections;  however, it may contain some exceptional curves with
multiplicity greater than one.  In this case we can apply
(\ref{smoothingtwo}) to decrease their total multiplicities.  (From another
perspective, the virtual dimension for
holomorphic sections of $X_2(f)$, in the homology class defined by the
section $[\mu]$ of divisors of multiplicity two supported along the exceptional
curve $E \subset X'$, is 
negative.  Hence for generic $J$ we do not expect the situation to
arise, although we do not need to prove that here.)  An easy
induction, pushing off multiple exceptional curves, yields a positive
symplectic divisor in $X'$ which is of the form $\bigcup E_i \cup D$,
where $D$ represents the class $p^* K_X$, with $p: X' \rightarrow X$
the blow-down map.  We still have to smooth $D$.  The last stage,
then, in the argument is to observe that the
conditions of (\ref{smoothinglemma}) are implied by the adjunction
formula. For given the symplectic surface $D = \bigcup a_i
\mathcal{C}_i$ we have that 

\begin{Eqn} \label{adjunctionequation}
K \cdot \mathcal{C}_i + (\mathcal{C}_i)^2 = 2g(\mathcal{C}_i) - 2 =
(1+a_i) (\mathcal{C}_i)^2 + \sum_{j \neq i} a_j (\mathcal{C}_i \cdot
\mathcal{C}_j).
\end{Eqn}

\noindent From these expressions, and the definition of a positive
symplectic divisor, if $\mathcal{C}_i ^2 \geq 0$ then $K \cdot (\mathcal{C}_i)
\geq 0$.  On the other hand, if we assume that both $(\mathcal{C}_i)^2
< 0$ and $K \cdot (\mathcal{C}_i) < 0$ then it follows that
$\mathcal{C}_i$ is a $(-1)$-sphere (and $X$ itself was not minimal).
Using (\ref{smoothingtwo}) we can again
separate out the $(-1)$-sphere components, to leave a positive
symplectic divisor which does satisfy the hypothesis of
(\ref{smoothinglemma}). 
Accordingly, we can smooth the components of $C_{\phi} \backslash \cup
E_i$ in a small neighbourhood of their union.  The resulting surface
remains disjoint 
from the exceptional sections $E_i$, and hence we can identify the
smooth symplectic
submanifold that results with a symplectic submanifold inside $X$, by
(\ref{symplecticforms}).  This represents $K_{X'} - \sum E_i = K_X$ in
homology, and the proof is complete.

\vspace{0.2cm}

\noindent To end, it may be worth 
mentioning that the hypothesis that the pencil has high degree
is indeed playing a definite role here in excluding bubbles.

\begin{Example} There is a genus two Lefschetz
pencil, with total space $X$, leading to a fibration $X'$ with mapping
class group word $(\delta_1 
\delta_2 \delta_3 \delta_4)^{10} = 1$ in standard generators.  (This
pencil is described in \cite{ivanhodge} for instance.)   The
four-manifold underlying the pencil is a  minimal complex surface of
general type on the Noether
line.  This surface is simply connected and has $b_+ > 1$ and hence
our basic index problem 
has solutions.  For the symplectic form on $X$ dual to one of the
genus two curves in the pencil, we have
$\omega^2 = 1, \, K_X \cdot \omega = 1$.  By the obstruction
computation, we have a section $\phi$ of a bundle of second symmetric
products over $\sss^2$; one of the two points on each fibre is
moreover that defined by the unique exceptional section $E$. 

\vspace{0.2cm}

\noindent Now if the other component of the cycle $C_{\phi}$ defined
by $\phi$ is a section of $f$ disjoint from $E$, then we have
represented $K_X$ by
a symplectic sphere, which violates the adjunction inequality.  It follows that
$[C_{\phi}] = [2E + F]$ for a fibre $[F] = [p^* K_X - E]$; indeed
for the underlying minimal surface, $K_X = [\omega]$ can indeed be
represented by a genus two curve, but the pencil of genus two curves
does not extend to a web.  The section of $X_2(f)$ is a cusp section,
and the curve $C_{\mathrm{bubble}} \subset X'$ pushes down to give a
surface in $X$, in the class $K_X$, which passes through the basepoint
of the pencil.

\vspace{0.2cm}

\noindent The point here is that $K_X = [\omega/2\pi]$ for this pencil of
low degree $k=1$.  It follows that there is a holomorphic representative for
$K_{X'} - [\Fibre]$ - the empty holomorphic curve - and this, stabilised
by the exceptional section, is trapping the section component of our
cusp curve.  The bubble defines a symplectic surface in the
four-manifold $X'$ which in fact projects to a smooth complex curve in
the class $K_X$. 
\end{Example}


\section{Appendix}

Our purpose here is to give a ``users' guide'', for non-specialists,
to the theory of divisors and line bundles on a nodal curve.  There is
a large literature on these topics to which we defer for careful
proofs.  One construction of the compactified Jacobian of a nodal 
curve, and of the relative Hilbert scheme (which gives the smooth
compactification of the family of symmetric products), involves
geometric invariant theory.
Standard (but inexhaustive) references include
\cite{Oda-Seshadri} for Jacobians of stable
curves, following on from work of Mumford \cite{Mumfordjacobians};
the relative Hilbert scheme is treated carefully in
\cite{moduliofshaves}.  A fine overview of Hilbert schemes on complex
surfaces is Nakajima's book \cite{Nakajima}.  Finally, a technical
survey covering all that we quote and more
is Kleiman's summary \cite{Kleiman}.  To fit more closely with these
references, in this Appendix we introduce $X_r (f)$ and $P_r (f)$
as pullbacks of the ``universal'' families over the moduli space of
curves;  the smooth structures one obtains this way are of course the
same as those induced by the restricted charts of
(\ref{restrictedcharts}).

\subsection{The projective bundle}

Let $\pi:Z \rightarrow B$ be
any holomorphic family of Riemann surfaces.  Then there is a unique vector
bundle $W \rightarrow B$ whose fibre over $b \in B$ is
\emph{canonically} identified with the space of holomorphic sections
$H^0 (K_{\pi^{-1}(b)})$.  For we can define a line bundle - the
dualising sheaf\footnote{The notation $\omega$ for dualising sheaves
  is as established as the identical notation for symplectic forms; we
  hope no confusion will arise.} - on $Z$ by 

$$\omega_{Z/B} \ = \ K_Z \otimes (\pi^* K_B)^{-1}$$

\noindent and then set $W = \pi_* \omega_{Z/B}$.  The bundle $W
\rightarrow B$ is called the \emph{relative dualising sheaf}.
 If we apply this construction to the universal
 curve $\pi_g: \mathcal{C}_g \rightarrow M_g$ then the bundle
 $(\pi_g)_* \omega$ 
 extends over the stable compactification $\mgbar$ \cite{moduliofcurves}.

\begin{Defn} \label{projbundle}
Let $f: X' \rightarrow \sss^2$ be a Lefschetz fibration inducing a map
$\phi_f: \sss^2 \rightarrow \mgbar$.  Let $V = \phi_f ^* ((\pi_g)_*
\omega) \otimes \mathcal{O}(-2)$ be the \emph{bundle of fibrewise canonical
  forms} associated to $f$.  This is a rank $g$ complex vector bundle
over $\sss^2$. 
\end{Defn}

\noindent The vector bundle $V \rightarrow \sss^2$ has the following
property.  A section of $V$ defines a cycle in $X'$
as follows: for each $b \in \sss^2$ we have a one-form on $f^{-1} (b)$
defined up to scale; the zeroes of this one-form are
well-defined and give a collection of $2g-2$ points, counted to
multiplicity, in $\pi^{-1}(b)$.  If the section of $V$ has no zeroes
we obtain a cycle in the class $PD[K_{X'}]$.  
We make a few remarks about this construction at the
critical values of $f$.  If $F \subset Z$ is a smooth complex curve in
a complex surface $Z$ (not necessarily compact) then there is an
adjunction formula for $K_F$:

\begin{Eqn} \label{Kfornodalcurves}
K_F \ = \ K_Z |_F \otimes \nu_{F/Z}
\end{Eqn}

\noindent where the last term denotes the normal bundle.  If $F$ is
now a nodal complex curve, it still defines an effective divisor and
hence line bundle $\mathcal{O}_Z(F)$ on $Z$, and we may formally
define a normal bundle $\nu_{F/Z} = \mathcal{O}_Z(F) |_F$ and hence a
canonical bundle by (\ref{Kfornodalcurves}).  The resulting locally
free sheaf is independent of the choice of $Z$ and of embedding $F
\subset Z$, as in the smooth case.  This defines an extension of the
dualising sheaf over $\mgbar$.  To give a geometric picture of the
elements of the fibre of $V$ over critical values of $f$ we can look
to the normalisation of the nodal fibre $\tilde{\Sigma}_0 \rightarrow
\Sigma_0$.  Recall that this is a naturally associated Riemann surface
in which the two sheets which meet at the node are separated.  If the
node does not separate $\Sigma_0$ then the normalisation has genus
smaller than that of a smooth deformation (connect sum) at the node,
and if the node separates then the normalisation is a disconnected
surface.  There are two
distinguished points $\alpha, \beta \in \tilde{\Sigma}_0$ given by the
preimages of the node.  

\begin{Prop} \label{residuemaps}
The elements of $H^0 (\Sigma_0, \omega_{\Sigma_0})$ can be identified
with the meromorphic sections of $K_{\tilde{\Sigma}_0}$ over
$\tilde{\Sigma}_0$ which are smooth away from $\alpha, \beta$ and have at
worst simple poles, with opposite residues, at each of $\alpha, \beta$.
\end{Prop}

\noindent This is clear from the explicit ``residue map'' \cite{BPV}; for a
local section of $K_Z \otimes \mathcal{O}_Z (\Sigma_0)$, say $h (du
\wedge dv) / f$ with $u,v$ local complex co-ordinates on $Z$, $h$ a
local section of $K_Z$ and $f$ a
local defining equation for $\Sigma_0$, we write 

$$(\mathrm{res}) \left( \frac{h \, du\wedge dv}{f} \right) \ = \
(\mathrm{norm})^* \left(  
\frac{h \, dv}{\cd f / \cd u} \right)$$

\noindent where the partial derivative is chosen not to vanish on any
open set in $\Sigma_0$.  Then $(\mathrm{res})$ identifies sections of
the canonical sheaf of $\Sigma_0$ with meromorphic forms on the
normalisation, and the poles arise from the zeroes of the derivative
$\cd f / \cd u$; these are at worst simple when $f$ is quadratic, as near a
node.  (Indeed the condition on the residues being opposite is forced
in this case, since the sum of residues must be trivial by Cauchy's theorem.)
For a smooth Riemann surface, the canonical bundle has degree $2g-2$
and any holomorphic section defines a distinguished set of $2g-2$
points (to multiplicity) via its zeroes.  The same is true for a nodal
Riemann surface;  a section $s \in H^0 (\Sigma_0, \omega_{\Sigma_0})$,
viewed as a meromorphic section on the normalisation, defines the
points

\begin{enumerate} 
\item which are zeroes of the meromorphic section on
  $\tilde{\Sigma}_0$ if it has poles at the points $\alpha, \beta$;
\item which are zeroes of the meromorphic section and the nodal point
  of $\Sigma_0$ if the section is actually holomorphic upstairs.
\end{enumerate}

\noindent This is again clear from the residue; if the residue form is
  smooth at the node, then the canonical section $h$ must vanish
  there.  It follows that in the case of a
  separating node, all the elements of $H^0 (K_{\tilde{\Sigma}_0})$
  give rise to tuples of points including the node; in algebraic
  geometry, the nodal points in reducible curves form base-points for
  the canonical linear system.  This explains why some arguments are
  simpler in the absence of reducible fibres.

\subsection{The Picard fibration}

There is a fibre bundle over $M_g$ with fibre the Picard torus of
degree $r$ line bundles on the associated Riemann surface.  (The
Picard variety of a complex curve identifies, after choosing an
origin, with the Jacobian of the curve.)  This also
extends to the stable compactification, but no longer as an orbifold.
For curves with a single node the situation is simpler than for
arbitrary stable curves.
The precise notion in algebraic geometry is the relative moduli
scheme for rank one torsion free coherent sheaves\footnote{Recall that
  any torsion 
free sheaf on a smooth curve is locally free, so the fibres do
reproduce the Picard varieties over the smooth locus $M_g$.} with fixed Euler
characteristic.  This is reduced
for stable 
curves (\cite{Oda-Seshadri}, Cor. 13.3);  when there is a unique node,
the compactified Jacobian can be described explicitly (below).  The
papers of Igusa \cite{Igusa}  gave an
early construction of compactified Jacobians using linear systems,
taking closures under suitable projective embeddings.  Mumford
\cite{Mumfordjacobians} showed
that the resulting spaces are indeed (components of) parameter spaces
for torsion free sheaves.

\begin{Defn}
Given a Lefschetz fibration $f:X' \rightarrow \sss^2$ write $P_r
(f)$ for the pullback by $\phi_f: \sss^2 \rightarrow \mgbar$ of the
relative moduli scheme of rank one torsion free sheaves of fixed Euler
characteristic $r-g+1$.  
\end{Defn}

\noindent Here is a description of this object near the singular fibres
of $f$.
A line bundle of degree $r$ on a nodal curve $\Sigma_0$ is given
by a line bundle of degree $r$ on the normalisation $\tilde{\Sigma}_0$
together with a 
choice of identification of the complex lines $L_{\alpha}$ and
$L_{\beta}$ over the preimages of the node $Q \in \Sigma_0$.  The
bundle upstairs is given by pullback: writing $(\mathrm{norm}):
\tilde{\Sigma}_0 \rightarrow \Sigma_0$ then take

$$\tilde{L} \ = \ (\mathrm{norm})^* L \otimes_{\mathcal{O}(\Sigma_0)}
\mathcal{O} (\tilde{\Sigma}_0) / \langle \mathrm{Torsion} \rangle;$$

\noindent dividing by torsion ensures the final sheaf is locally free.
This
yields a non-compact space, which is a $\cc^*$ fibre bundle over the
Picard torus of the normalisation.  The natural compactification of this
to a $\pp^1$-bundle, allowing the degenerate gluing maps of the two
complex lines by the $0$ and $\infty$ multiplications, is no longer a
moduli space for a natural class of objects.  However, if we
parametrise rank one torsion free sheaves on $\Sigma_0$, then we find
that a quotient of this $\pp^1$ bundle is the required moduli space.
Precisely, we glue together the $0$-section and $\infty$-section of
the $\pp^1$ bundle over an automorphism of the base torus which is a
translation by $\mathcal{O}(\alpha - \beta)$ in the group action of
degree zero line bundles on the Picard.  This is derived in
(\cite{Oda-Seshadri}, Example p.83) and also (\cite{Igusa},
Supplement p.187).

\vspace{0.2cm}

\noindent One can see the degeneration of the Picard / Jacobian
fibration explicitly in terms of periods.  Suppose we have a smooth
complex surface $f: \mathcal{X} \rightarrow D$ with a nodal fibre
$\Sigma_0$ over $0$, with normalisation $\tilde{\Sigma}_0$ as usual.
Fix a basis of loops for 
$H_{1}(\Sigma_{t})$ in the obvious way:
$\gamma_{3}\dots, \gamma_{2g}$ correspond to a standard basis on
$\tilde{\Sigma}_0$,
$\gamma_{1}$ is the vanishing cycle, $\gamma_{2}$ a loop with 
$\gamma_{2} \cdot \gamma_{1}=1$. Note that because of the topological
monodromy we  
can only define $\gamma_{2}$ consistently on a covering or cut plane. Now let 
$\omega_{\alpha}$ be a basis for the ``$1$-forms'' on the singular 
fibre, which we can obtain by the residue map from holomorphic $1$-forms on 
$\mathcal{X}$. We can suppose that the $\omega_{\alpha}$, for $\alpha>1$, are 
holomorphic on the normalisation while $\omega_{1}$ has residue $1$.
Moreover, by the discussion of the residue map above, we can regard
$\omega_{\alpha}$ as being defined on all the 
curves $\Sigma_{t}$.
We consider the periods

$$   \int_{\gamma_{i}} \omega_{\alpha}. $$

\noindent It is not hard to see that these are all holomorphic
functions of $t$ except that

$$  \int_{\gamma_{2}} \omega_{1} = \log t \int_{\gamma_{1}} \omega_{1}
+ \mathrm{holomorphic}.  $$

\noindent Moreover at $t=0$ the integrals
$\int_{\gamma_{i}}\omega_{\alpha}$ for 
$i\geq 3$ and $\alpha\geq 2$ give the periods of $\tilde{\Sigma}_0$, while
$\int_{\gamma_{1}} \omega_{1} = 2\pi $ and $\int_{\gamma_{1}}
 \omega_{\alpha} =0$ for $\alpha\geq 2$.

\vspace{0.2cm}

\noindent  This gives an explicit description of the total space of
$\mathcal{P}_r(f)$ (induced from $\mathcal{X} \rightarrow D$) minus the  
normal crossing divisor, as a quotient of $\cc^{g}\times D$. Then one can 
construct the compactification of this explicitly in local co-ordinates. The 
basic model is to consider the quotient of $\cc\times D$ by the 
equivalence relation 

$$ (z, t) \sim (z + 2 \pi (n_{1} + n_{2} i \log t), t), $$

\noindent if $t\neq 0$ and 

$$(z,0)\sim (z + 2\pi n_{1},0).$$

\noindent Consider the regions $0< \Im(z)<-\Re(\log t)$ and $\Re(\log
t)< \Im (z)<0$, where $\Re$ and $\Im$ denote real and imaginary parts
respectively. 
We map the first region to $\cc^{2}$ by $(z,t)\mapsto
(e^{iz},te^{-iz})$ and the  second by $(z,t)\mapsto (te^{iz}, e^{-iz})$. 
These induce a map from a neighbourhood of the end of the quotient to 
$\cc^{2}\backslash \{0\}$, and we compactify by adding $0$.

\begin{Prop}
For a Lefschetz fibration $f: X' \rightarrow \sss^2$ with irreducible
fibres, the Picard fibration $P_r (f)$ has smooth symplectic total
space.  The critical points of the natural map $P_r (f) \rightarrow
\sss^2$ are precisely the normal crossings divisors described above in
the fibres $\Pic_r (\Sigma_{f(p_i)})$, for $\{ f(p_1), \ldots, f(p_n) \}$ the
critical values of $f$.
\end{Prop}

\noindent The global smoothness will follow from the discussion for
symmetric products below, or can be proven directly in an analogous fashion.
For irreducible curves with one node, the natural tensor
product action of
degree zero line bundles on the compactified Picard has a single
orbit.  Under the map given by $\otimes \omega_{\Sigma_0}$ from
$\Pic_0$ to $\Pic_{2g-2}$ the point corresponding to
$(\mathcal{O}(\tilde{\Sigma}_0), 1 \in \cc^*)$ is mapped to the
canonical sheaf of the nodal curve.  The locus
of critical values for the natural projection to $\sss^2$ corresponds to
torsion free non-locally free sheaves.  The canonical sheaf of a nodal
curve \emph{is} locally free; hence the natural section of $P_r (f)
\rightarrow \sss^2$ defined by taking a point $t \in \sss^2$ to the
canonical sheaf $K_{f^{-1} (t)}$ is well-defined and smooth.

\subsection{The relative Hilbert scheme}

The symmetric product of a smooth curve is equal to the Hilbert
scheme, parametrising fixed length quotients of the structure sheaf.  For
background on Hilbert schemes and Quot schemes see
\cite{moduliofshaves}, where relative 
Hilbert schemes are shown to exist as projective schemes in great
generality.  (From the algebraic perspective, after fixing appropriate
discrete data, symmetric products
are moduli spaces of structure sheaves whilst Hilbert schemes are
moduli spaces of ideal sheaves.) 
For a family of curves with isolated nodal members, this gives a
compactification of the fibre bundle of symmetric products.  The
important points for us will be that the total space is smooth and
there is still a  well-defined and
global Abel-Jacobi map.   A careful treatment of the latter assertion
can be found in 
Altman and Kleiman (\cite{Altman-Kleiman}, Section 8), and we shall
provide a proof of the former (\ref{globalsmoothness}).  Some additional
information on these spaces is given in \cite{Sequel}.

\begin{Defn}
Given a Lefschetz fibration $f: X' \rightarrow \sss^2$ with
irreducible fibres, write $F: X_r (f) \rightarrow \sss^2$ for the
total space of the relative Hilbert scheme.  This is the
$\phi_f$-pullback of the scheme $\Hilb^{[r]}(\mathcal{C}_g / M_g)$
which parameterises length $r$ subschemes of the fibres of the
universal curve.
\end{Defn}

\noindent Given $f: X' \rightarrow \sss^2$ smooth over a locus
$\sss^2_*$, construct $P^*_r(f)$ over $\sss^2_*$.  By the previous remarks this
extends naturally to the entire sphere.   As for the smooth fibres, we
can identify a family of projective spaces over the 
Picard torus of the nodal curve.  Let $P \in \Pic_r (\Sigma_0)$ be a point
of the smooth locus, arising from a line bundle $L_P \rightarrow
\tilde{\Sigma}_0$ on the normalisation together with a gluing
parameter $\lambda \in \cc^*$ to identify the $\alpha$ and $\beta$
fibres of $L_P$.  Suppose for definiteness that the normalisation is
connected of genus $g-1$.  Then $L_P$, having degree $2g-2$,
generically has a space of holomorphic sections of dimension 

$$(2g-2) - (g-1) + 1 \ = \ g.$$

\noindent We take the $\lambda$-hyperplane in the space of all holomorphic
sections; that is, restrict to sections whose values at $\alpha$ and
$\beta$ are transformed by the gluing $\lambda: (L_P)_{\alpha}
\rightarrow (L_P)_{\beta}$.  Hence for a generic point $P$
and line $L_P$ the projective space of sections is a copy of
$\pp^{g-2}$.  A similar analysis applies along the normal crossing
divisor.  This changes, however, at a unique point in the smooth
locus of $\Pic_r (\Sigma_0)$, corresponding to the canonical line of the
normalisation.  Here we see a projective space of dimension $g-1$ as
the space of sections.  Note this description of the space of sections
$H^0 (\omega_{\Sigma_0})$ differs from that given above in terms of
meromorphic differentials:  $(\norm)^* \omega_{\Sigma_0} / \langle
\mathrm{Tors} \rangle$ is \emph{not} the canonical sheaf of the
normalisation, it has the wrong degree.  

\vspace{0.2cm}

\noindent Following work of Nakajima, we can give a more down-to-earth
description of the
relative Hilbert scheme for the local model $\pi: \cc^2 \rightarrow \cc$
taking $(z,w) \mapsto zw$.  For $\Hilb^{[r]}$ the fibre over $t$ is
the set of ideal sheaves 
in the local ring $\cc[z,w]/ \langle zw-t \rangle$ whose quotient is
of length $r$.  Equivalently we want the ideals in $\cc[z,w]$ which
contain $\langle zw-t \rangle$ and which have quotient length $r$.
The Hilbert scheme of the complex plane has an elementary description
\cite{Nakajima}: 

$$\Hilb^{[r]} (\cc^2) \ = \ \big \{ (B,B',v) \in M_r (\cc) \times M_r
(\cc) \times \cc^r \ | \ [B,B'] = 0, \ (*) \big \} \big / GL_r (\cc)$$

\noindent where $(*)$ is a stability condition for a geometric
invariant theory construction:  no subspace $S \subset \cc^r$
invariant under each of $B$ and $B'$ can contain the vector $v$.  The
gauge group $GL_r (\cc)$ acts by simultaneous conjugation on $B$ and
$B'$ and by left multiplication on $v$.
Given an ideal $\mathcal{I}$, the commuting matrices $B$ and $B'$ represent
multiplication by $z$ and $w$ respectively on the $r$-dimensional
vector space $\cc[z,w]/\mathcal{I}$, and $v$ arises as the image of $1
\in \cc[z,w]$.  Conversely, given a triple as above, we define a map
$\phi: \cc[z,w] \rightarrow \cc^r$ by $f \mapsto f(B, B') v$; the stability
condition ensures that $\phi$ is surjective and then the kernel defines an
ideal $\mathcal{I}_{\phi}$.  This gives a
set-theoretic description of 
the Hilbert scheme for length $r$ quotients, and \cite{Nakajima} shows
that this is a holomorphic isomorphism.  Again by stability, the ideal
associated to $(B, B', v)$ is just $\{ f \in \cc[z,w] \ | \ f(B,B') =
0 \}$.  This contains $zw-t$, for given $t \in \cc$, just when
$BB'=tI_r$.  In particular, this gives a fairly explicit description
of the singular fibre $\Hilb^{[r]} (\{zw=0\})$:

$$\Hilb^{[r]} (\pi^{-1}(0)) \ = \ \big \{ (B, B', v) \ | \ [B,B'] = 0, \ BB'
= 0, \ (*) \big \} / GL_r (\cc).$$

\noindent From this perspective, a
general point of the symmetric product of this fibre $\pi^{-1}(0)$ can
be given by a pair of diagonal matrices $diag(\lambda_1, \ldots,
\lambda_r)$ and $diag(\mu_1, \ldots, \mu_r)$ with $\lambda_i \mu_i =
0$ for each $i$;  a general point of the Hilbert scheme can be given
by a suitable pair of upper triangular matrices $(B,B')$ whose
diagonal entries satisfy 
the same condition.  The Hilbert-Chow morphism from $\Hilb^{[r]}
\rightarrow \Sym^r$ takes $(B,B')$ to the set of pairs $(\lambda_i,
\mu_i)$ viewed as points of $\cc^2$ lying on $\pi^{-1}(0)$.  Where the
eigenvalues of $B$ and $B'$ all remain distinct, this map is an
isomorphism; one can use this explicit form of the map to prove
(\ref{stillcodim2}).  The Hilbert scheme is stratified by the
number of supports $(\lambda, \mu)$ lying at the node of $\{ zw=0 \}$,
and as points collapse into the node configurations of off-diagonal matrix
entries become (projective) co-ordinates over the stratum.

\begin{Example}
Suppose $r=2$.  If either $B$ or $B'$ has distinct eigenvalues then we
can simultaneously diagonalise both with $B = diag(\lambda_1,
\lambda_2)$ and $B' = diag(\mu, \mu')$ and with $(\lambda_1, \mu_1)
\neq (\lambda_2, \mu_2)$.  This point of $\Hilb^{[2]} (\cc^2)$
belongs to the relative Hilbert scheme of $\pi$ precisely when
$\lambda_1 \mu_1 = \lambda_2 \mu_2$.  If, however, $B$ and $B'$ have
only one eigenvalue each, then they cannot be simultaneously
diagonalised by the stability condition $(*)$ (otherwise take $S =
\langle v \rangle$).  Hence we can write 

$$B = \left( \begin{array}{cc} \lambda & \alpha \\ 0 & \lambda
  \end{array} \right); \qquad B' = \left( \begin{array}{cc} \mu &
    \beta \\ 0 & \mu \end{array} \right)$$

\noindent for $(\alpha, \beta) \in \cc^2 \backslash \{0\}$.  The
associated ideal is generated by 

$$\big \langle \beta(z-\lambda) - \alpha(w-\mu), \ (z-\lambda)^2, \
(w-\mu)^2 \big \rangle.$$ 

\noindent This
subscheme is supported on the fibre $\pi^{-1}(0)$ when $\lambda=0$ or
$\mu=0$.  However, if $\lambda = 0$ and $\mu \neq 0$ then the ideal
does \emph{not} contain $\langle zw \rangle$, and the corresponding
quotient is not a point of the relative Hilbert scheme.  (Thus the
relative Hilbert scheme is not just the preimage of the relative
symmetric product under the Hilbert-Chow morphism.)  We will interpret
the parameters $\alpha$ and $\beta$ geometrically below.
\end{Example}

\noindent From this
description, one can prove the smoothness of all the relative Hilbert
schemes $X_r(f)$, which in addition implies the smoothness of the
fibre product spaces $Y_{\pi, \aleph}$.  The details of this are
carried out in \cite{Sequel}, but do not illuminate how the strata of
the relative 
Hilbert scheme fit together..  Instead, here is a more direct and
geometric argument for smoothness of $X_{2g-2}(f)$ - which can be
adapted to the general case - based upon our ``normal crossing''
picture~:

\begin{Prop} \label{globalsmoothness}
Let $f: \mathcal{X} \rightarrow D$ be a smooth complex surface which fibres
over the disc with an irreducible nodal fibre $\Sigma_0$ over
$0$, as usual.  The
total space of the relative Hilbert scheme $\mathcal{X}_r(f) \rightarrow D$ is
smooth.
\end{Prop}

\begin{Pf}
From the discussion above it follows that the total space has normal
crossing singularities in the fibre over $0 \in D$ whilst all the
other fibres are given by the (smooth) symmetric products of the
fibres of $f$.  The general theory of deformations of spaces with
normal crossings was carefully described by Friedman in
\cite{Friedman:NC}.  Be given a flat proper map $\mathcal{Y} \rightarrow
D$ over the disc, with central fibre $\mathcal{Y}_0 = Y$ having
normal crossing singularities.  There is a Kodaira-Spencer map 

$$\theta: (TD)_0 \ \rightarrow \ \mathbf{Ext}^1 _{\mathcal{O}_Y}
(\Omega^1 _Y, \mathcal{O}_Y)$$

\noindent to the global Ext-group, which is formally the tangent space
to the deformation space of $Y$.  The class $\theta( \partial /
\partial t)$ defines an extension which is just the conormal sequence

$$0 \rightarrow \mathcal{I}_Y / \mathcal{I}_Y ^2 \ \rightarrow \ (\Omega^1
_{\mathcal{Y}})|_Y \ \rightarrow \ \Omega^1 _Y \rightarrow 0.$$

\noindent Since $Y$ appears as a fibre in a flat family, the first
term above is $\mathcal{O}_Y$.  Now Friedman's theorem asserts
that $\mathcal{Y}$ is smooth at a point $y \in Y$ precisely where
$\theta (\partial / \partial t)$ generates the local
group $\mathcal{E}xt^1 _{\mathcal{O}_Y} (\Omega^1 _Y,
\mathcal{O}_Y)$.  (Note that the conormal sequence also defines an
element of $H^0 (\mathcal{E}xt)$.)  In our situation, this local
Ext-group is exactly 
$\mathcal{O}_R$, where $R$ denotes the singular locus of $Y$, and
where the sheaf $\mathcal{O}_R$ denotes 
the normal bundle of $Y$ in $\mathcal{Y}$ restricted to $R$.

\vspace{0.2cm}

\noindent More concretely, the germ of the family
$\mathcal{X}_r(f) \rightarrow D$ over $0$ 
defines a germ of a map $D \rightarrow \mathrm{Def}(\mathcal{X}_0) =
\mathbf{Ext}$ 
into the deformation space.  The derivative of this
deformation at $0$ defines a section $\epsilon$ over $\mathcal{X}_0$
of the line bundle $N_1 ^* \otimes N_2 ^*$, where the $N_i$ are the
normal bundles to the two branches of the normal crossing divisor.
Friedman's theorem \cite{Friedman:NC} amounts to saying that 
it suffices to check that the section $\epsilon$ is nowhere zero.

\vspace{0.2cm}

\noindent In our case, the normal bundles $N_1$ and $N_2$ can be
identified canonically with the tangent lines
$(T\tilde{\Sigma}_0)_{\alpha}$ and $(T \tilde{\Sigma}_0)_{\beta}$ to
the normalisation of the fibre of $f$ at the marked points.  (A
similar statement for curves is presented, in down-to-earth fashion, in
\cite{moduliofcurves}.)  It
follows that $\epsilon$ gives a constant section of a trivial line
bundle, and we need to see that we obtain a non-zero element of the
line.  Conversely, if the section is trivial then the total space will
be singular along the entire normal crossing divisor, so it suffices
to check smoothness at a single point.  Choose $2g-3$ local sections
of $\mathcal{X} \rightarrow D$;  using the semigroup structure on
$\oplus_d \Sym^d$ we obtain a map $\mathcal{X}_1(f) \rightarrow
\mathcal{X}_{2g-2}(f)$ 
which adds twice the first co-ordinate to the given sections.
However, $\mathcal{X}_1(f) = \mathcal{X}$ is exactly the complex surface we
start with.  It follows that we have a factorisation

\begin{Eqn} \label{factorise}
\mathcal{X} \stackrel{\phi}{\longrightarrow} \mathcal{X}_{2g-2}(f)
\rightarrow D. 
\end{Eqn}

\noindent Let $x,y$ be co-ordinates near the node in $\mathcal{X}$ so
that locally the map $f$ is given by $(x,y) \mapsto xy$;  let $z_1,
\ldots, z_n$ be co-ordinates near a point of the normal crossing
divisor in $\mathcal{X}_r(f)$, so the projection and section are given by

$$(z_1, \ldots, z_n) \mapsto z_1 z_2; \qquad z_1 z_2 = \epsilon(z_3,
\ldots, z_n)$$

\noindent respectively.  Then (\ref{factorise}) implies that locally
we can write the product

$$z_1 (x,y) z_2(x,y) \ = \ F(xy)$$

\noindent where $\frac{\partial z_1}{\partial x}$ and $\frac{\partial
  z_2}{\partial y}$ are both non-zero.  It follows that $F'(0) \neq
  0$, and hence that $\mathcal{X}_r(f)$ is smooth at points in the normal
  crossing divisor in the image of $\phi$, and hence smooth globally.
\end{Pf}

\noindent The existence of a relatively ample holomorphic line bundle on the
  total space of $\mathcal{X}_r(f)$, extending that induced over the punctured
  disc by a given collection of sections, follows from
  (\cite{moduliofshaves}, 2.2.5).  Roughly, the line bundle we want is
  given by a twist of the determinant line bundle of an
  appropriate 
  ``universal'' line.  The geometric invariant theory construction of
  Hilbert schemes proceeds by embedding the required set of quotients
  of a sheaf $\mathcal{H}$ over a space $X$ 
  into a Grassmannian of subspaces of $H^0 (\mathcal{H} \otimes
  L_X ^m)$, where we
  have twisted $\mathcal{H}$ by some ample sheaf $L_X$ on
  the space $X$.  The determinant of the universal line is ample
  since it is the
  pullback of the $\mathcal{O}(1)$ bundle on some $\pp^N$ under the
  Pl{\"u}cker (determinant) embedding of the Grassmannian.  In the
  relative case, we twist by a bundle which is ample on the fibres of
  the given map of schemes $X \rightarrow B$, which certainly holds
  for the line bundle on a Lefschetz fibration defined by the
  exceptional sections.  This completes our treatment of 
  the material underlying Theorem (\ref{itallextends}).

\begin{Example}
For a smooth genus two curve, the second symmetric product maps to the
Picard by blowing down a single $\cc \pp^1$ which is the fibre over
the canonical line bundle.  
Given a genus two curve $\Sigma_0$ with a single non-separating node,
the second 
symmetric product is singular along a copy of $\Sigma_0$ parametrising
all pairs $(p, \mathrm{node})$ with $p \in \Sigma_0$.  The Hilbert
scheme for two points on $\Sigma_0$ is given by blowing up the second
symmetric product at the point $(\mathrm{node}, \mathrm{node})$ which
creates an exceptional divisor $E \cong \cc \pp^1$.  (In the
description above by ideals, the point $[\alpha: \beta]$ defines a
co-ordinate on this copy of $\pp^1$.)  The
singular locus of $\Hilb^{[2]} (\Sigma_0)$ is a copy of the normalisation
$\tilde{\Sigma}_0$; the generic point on the exceptional divisor $E$
is smooth.  The two singular points on $E$ correspond to subschemes
with the two
points lying at the node, with an infinitesimal deformation (showing
the direction in which they collided) also being tangent to the node.
Locally, the blow-up serves to give a crepant resolution of the
(globally singular) fibrewise second symmetric product of $\pi:\cc^2
\rightarrow \cc$ taking $(z,w) \mapsto zw$.
\end{Example}

\noindent We end the Appendix with a proof of the technical result
(\ref{bubblingfacts}) asserting that bubbles in the singular fibres of
$F$ realise no more homology classes than bubbles in the smooth fibres.

\begin{Lem} \label{proofforbubbles}
Let $\theta: \sss^2 \rightarrow \Hilb^{[r]} (\Sigma_0)$ be a non-constant
holomorphic map from the two-sphere to the singular fibre of $F: X_r(f)
\rightarrow D$ over $0 \in D$ which is a bubble component of some cusp
curve.  Then $[\im(\theta)] \equiv Nh$; the homology class of the image
is equal to a multiple of the class $h$ defined by a projective line
inside $\pp(H^0 (K_{\Sigma_t})) \subset F^{-1}(t)$.
\end{Lem}

\begin{Pf}
Compose $\theta$ with the map $\tau$ from
$X_r(f)$ to $P_r(f)$; we get a map into $\Pic_r(\Sigma_0)$. This lifts
to the normalisation of 
$\Pic$ which we have seen fibres over $\Pic(\tilde{\Sigma}_0)$. Hence
the composite 
$\sss^{2}\rightarrow \Pic(\tilde{\Sigma}_0)$ is constant. The only
possibilities are that 
either the map $u=\tau \circ \theta: \sss^{2}\rightarrow \Pic_r (\Sigma_0)$ is
constant or it is maps non-trivially 
onto one of the $\sss^{2}$ fibres of $\Pic_r  (\Sigma_0)  \rightarrow
\Pic(\tilde{\Sigma}_0)$.  In the first 
case the map $u$ is clearly homotopic to the standard sphere so we have to 
rule out the second case. In this case $\theta$ must meet the singular
locus of  
$\Hilb^{[r]}(\Sigma_0)$ at least twice. We claim that this means that
$\theta$ cannot arise as the   
bubble component of a sequence of sections $\phi_{i}$ of the fibration
$X_r(f)$. If it did 
there would be pairs of disjoint discs $A_{i}, B_{i}\subset D$ and discs 
$A,B\subset \sss^{2}$ centred on points $a,b$ such that
$\theta(a), \theta(b) \in \{\mathrm{Singular \ locus}\}$, and such that after 
rescaling $\phi_{i}\vert_{A_{i}}$ converges in $C^{\infty}$ to
$\theta |_{A}$ (respectively for the $B$'s). 
Now consider the situation in local co-ordinates 
around $\theta(a)$.  From our identification of the singularities as normal
crossing divisors in the fibre, the projection map is given by
$$ (z_{1},z_{2}, \dots, z_{n})\mapsto z_{1} z_{2}. $$
The map $\theta$ must map locally into one of the branches $z_{1}=0,
z_{2}=0$ of  
the singular fibre, say into the first branch. This means that $\phi_{i}$ has 
intersection number $\geq 1$ with the hyperplane $z_{1}=0$ so there is a point 
$\alpha_{i}$ in $A_{i}$ such that $\phi_{i}(\alpha_{i})$ lies in the fibre 
of $X_r(f)$ over $0$. Similarly there is another point in $B_{i}$ with
the same  
property. Hence $\phi_{i}$ meets this fibre at least twice, and cannot
be a section, which gives the contradiction we require.
\end{Pf}


\bibliographystyle{alpha}
\bibliography{main}

\end{document}